\magnification=1200
\input amssym.def
\input amssym.tex

\hsize=13.5truecm
\baselineskip=16truept
\baselineskip=16truept
\font\secbf=cmb10 scaled 1200
\font\eightrm=cmr8
\font\sixrm=cmr6

\font\eighti=cmmi8

\font\sixi=cmmi6
\skewchar\eighti='177 \skewchar\sixi='177

\font\eightsy=cmsy8
\font\sixsy=cmsy6
\skewchar\eightsy='60 \skewchar\sixsy='60

\font\eightit=cmti8

\font\eightbf=cmbx8
\font\sixbf=cmbx6

\let\sc=\tensc

\font\eightsc=cmcsc10 scaled 800
\font\secbf=cmb10 scaled 1200
\font\subsecfont=cmb10 scaled \magstephalf
\font\amb=cmmib10

\font\ambi=cmmib10 scaled 700

\newfam\mbfam 

\textfont\mbfam\amb \scriptfont\mbfam\ambi


\def\aa{\def\rm{\fam0\eightrm}%
  \textfont0=\eightrm \scriptfont0=\sixrm \scriptscriptfont0=\fiverm
  \textfont1=\eighti \scriptfont1=\sixi \scriptscriptfont1=\fivei
  \textfont2=\eightsy \scriptfont2=\sixsy \scriptscriptfont2=\fivesy
  \textfont3=\tenex \scriptfont3=\tenex \scriptscriptfont3=\tenex
  \def\sc{\eightsc}
  \def\it{\fam\itfam\eightit}%
  \textfont\itfam=\eightit
  \def\bf{\fam\bffam\eightbf}%
  \textfont\bffam=\eightbf \scriptfont\bffam=\sixbf
   \scriptscriptfont\bffam=\fivebf
  \normalbaselineskip=9.7pt
  \setbox\strutbox=\hbox{\vrule height7pt depth2.6pt width0pt}%
  \normalbaselines\rm}

\def\Proof{\vskip12pt\noindent{\bf Proof.} }

\def\Def#1{\vskip12pt\noindent{\bf Definition #1}}
\def\Remark#1{\vskip12pt\noindent{\bf Remark #1}}
\def\Conclusion#1{\vskip12pt\noindent{\bf Conclusion #1}}

\def\m@th{\mathsurround=0pt}

\def\cc#1{\hbox to .89\hsize{$\displaystyle\hfil{#1}\hfil$}\cr}
\def\lc#1{\hbox to .89\hsize{$\displaystyle{#1}\hfill$}\cr}
\def\rc#1{\hbox to .89\hsize{$\displaystyle\hfill{#1}$}\cr}

\def\eqal#1{\null\,\vcenter{\openup\jot\m@th
  \ialign{\strut\hfil$\displaystyle{##}$&&$\displaystyle{{}##}$\hfil
      \crcr#1\crcr}}\,}

\def\section#1{\vskip 22pt plus6pt minus2pt\penalty-400
        {{\secbf
        \noindent#1\rightskip=0pt plus 1fill\par}}
        \par\vskip 12pt plus5pt minus 2pt
        \penalty 1000}

\def\subsection#1{\vskip 20pt plus6pt minus2pt\penalty-400
        {{\subsecfont
        \noindent#1\rightskip=0pt plus 1fill\par}}
        \par\vskip 8pt plus5pt minus 2pt
        \penalty 1000}

\def\subsubsection#1{\vskip 18pt plus6pt minus2pt\penalty-400
        {{\subsecfont
        \noindent#1}}
        \par\vskip 7pt plus5pt minus 2pt
        \penalty 1000}

\def\center#1{{\begingroup \leftskip=0pt plus 1fil\rightskip=\leftskip
\parfillskip=0pt \spaceskip=.3333em \xspaceskip=.5em \pretolerance 9999
\tolerance 9999 \parindent 0pt \hyphenpenalty 9999 \exhyphenpenalty 9999
\par #1\par\endgroup}}

\def\\{\hfill\break}

\def\kwadrat{\null\ \hfill\null\ \hfill$\square$}
\def\mida#1{{{\null\kern-4.2pt\left\bracevert\vbox to 6pt{}\!\hbox{$#1$}\!\right\bracevert\!\!}}}
\def\midy#1{{{\null\kern-4.2pt\left\bracevert\!\!\hbox{$\scriptstyle{#1}$}\!\!\right\bracevert\!\!}}}

\def\diagint{{\raise1.5pt\hbox{$\scriptscriptstyle\diagup$}\hskip-8.7pt\intop}}

\def\divv{{\rm div}\,}

\def\supp{{\rm supp}\,}

\def\today{${\scriptscriptstyle\number\day-\number\month-\number\year}$}
\footline={{\hfil\rm\the\pageno\hfil${\scriptscriptstyle\rm\jobname}$\ \ \today}}

\def\D{{\Bbb D}}

\def\R{{\Bbb R}}
\def\N{{\Bbb N}}
\def\T{{\Bbb T}}
\def\Z{{\Bbb Z}}

\def\esssup{{\rm esssup\,}}

\center{\secbf A priori estimate for axially symmetric solutions to the 
Navier-Stokes equations near the axis of symmetry}
\vskip1.5cm

\centerline{\bf Wojciech M. Zaj\c aczkowski}

\vskip1cm
\noindent
Institute of Mathematics, Polish Academy of Sciences,\\
\'Sniadeckich 8, 00-956 Warsaw, Poland\\
E-mail:wz@impan.pl;\\
Institute of Mathematics and Cryptology, Cybernetics Faculty,\\
Military University of Technology, Kaliskiego 2,\\
00-908 Warsaw, Poland
\vskip0.8cm


\vskip1.5cm

\noindent
{\bf Abstract.} 
We consider the axially symmetric solutions with large swirl to the 
Navier-Stokes equations. Let $v_r$, $v_\varphi$, $v_z$ be the cylindrical 
coordinates of velocity and $\chi=v_{r,z}-v_{z,r}$ be a component of 
vorticity. Let $V_2^k(\Omega^T)$ be a space with the finite norm
$$
\|w\|_{V_2^k(\Omega^T)}=\|w\|_{L_\infty(0,T;H^k(\Omega))}+
\|\nabla w\|_{L_2(0,T;H^k(\Omega))},
$$
$k=0,1$. We proved the a priori estimate for 
$\|\chi\|_{V_2^0(\Omega_\zeta^T)}$, where 
$\Omega_\zeta=\{x\in\Omega:\ r\le2r_0\}$ and $r_0$ is so small that the 
swirl $u=rv_\varphi$ satisfies
$$
\|u\|_{C(0,T;C^{1/2}(\Omega_\zeta))}\le\root{4}\of{5\over4}\nu,
\leqno(1)
$$
where $\nu$ is the viscosity and $u$ vanishes on the axis of symmetry. (1) is 
proved under the a priori estimate that $v\in W_2^{2,1}(\Omega^T)$ (see [Z5]).

\noindent
Next, having that $\chi\in V_2^0(\Omega^T)$ (proof of the esimate for 
$\|\chi\|_{V_2^0((\Omega\setminus\Omega_\zeta)^T)}$ is not diffficult) the 
elliptic problem $v_{r,z}-v_{z,r}=\chi$, $v_{r,r}+v_{z,z}+{v_r\over r}=0$ 
implies that $v'=(v_r,v_z)\in V_2^1(\Omega^T)$, so $v'\in L_{10}(\Omega^T)$. 
The above construction implies the following step by step procedure. Let 
$v(0)\in H^1(\Omega)$ then $v\in W_2^{2,1}(\Omega^T)$. Next 
$\chi\in V_2^0(\Omega^T)$ implies that $v(T)\in H^1(\Omega)$, where 
$v_\varphi(T)\in H^1(\Omega)$ is shown separately. This will be a topic 
of the next paper (see [Z4]).

\section{1. Introduction}

The uniqueness and regularity of weak solutions to the Navier-Stokes equations 
are important open problems in the incompressible hydrodynamics. Since the 
problems are very difficult in a general case there are examined special 
solutions to the Navier-Stokes equations such as: 
two\--dimensional, axially symmetric, helicoidal. The problems were solved for 
two-dimensional solutions in [L3] and in the case of helicoidal solutions in 
[MTL]. Moreover, the existence of global regular solutions which are close 
to the two-dimensional solutions was proved in [RZ, NZ, Z6, Z7]. So far the 
existence of regular axially symmetric solutions (see Definition 1.1) is 
still an open problem. However, the existence of regular axially symmetric 
solutions with vanishing swirl was proved long time ago in [L2, UY] the case 
with nonvanishing swirl is still open. Moreover, the existence of solutions 
which remain close to solutions from [L2, UY] was shown in [Z1, Z2, Z3, Z8]. 
In this paper and in [Z4] the existence of global regular axially symmetric 
solutions with nonvanishing swirl is proved in a periodic cylinder and 
with the slip boundary conditions on the lateral part of its boundary. The 
aim of this paper is to establish an a priori estimate for regular solutions 
in a neighbourhood of the axis of symmetry. In [Z4] a similar estimate is 
proved in a neighbourhood located at a positive distance from the axis of 
symmetry. Then in [Z4] by an appropriate partition of unity the estimate is 
derived in a whole cylinder. Finally in [Z4] by the 
Leray-Schauder fixed point theorem the existence of regular global axially 
symmetric solutions is established in a~whole cylinder.

\noindent
In this paper we consider the axially symmetric solutions (see 
Definition 1.1) to the following problem
$$\eqal{
&v_t+v\cdot\nabla v-\nu\Delta v+\nabla p=0\quad &{\rm in}\ \ 
\Omega^T=\Omega\times(0,T),\cr
&\divv v=0\quad &{\rm in}\ \ \Omega^T,\cr
&v\cdot\bar n=0\quad &{\rm on}\ \ S^T,\cr
&\bar n\cdot\D(v)\cdot\bar\tau_\alpha=0,\ \ \alpha=1,2,\quad &{\rm on}\ \ 
S_1^T=S_1\times(0,T),\cr
&{\rm periodic\ boundary\ conditions}\quad &{\rm on}\ \ S_2^T=S_2\times(0,T),
\cr
&v|_{t=0}=v_0\quad &{\rm in}\ \ \Omega,\cr}
\leqno(1.1)
$$
where $\Omega$ is an axially symmetric cylinder with boundary $S=S_1\cup S_2$, 
$x=(x_1,x_2,x_3)$ is the Cartesian system of coordinates in $\R^3$ such 
that $x_3$-axis is the axis of the cylinder $\Omega$. By $S_1$ we denote the 
part of the boundary of the cylinder parallel to the $x_3$-axis and $S_2$ 
is perpendicular to it.
Next, $v=v(x,t)=(v_1(x,t),v_2(x,t),v_3(x,t))\in\R^3$ is the velocity of the 
considered fluid, $p=p(x,t)\in\R$ the pressure, $\nu>0$ is the viscosity 
coefficient, $\D(v)=\nabla v+\nabla v^T$ is the double symmetric part of 
$\nabla v$, $\bar n$ is the unit outward normal vector to $S_1$ and 
$\bar\tau_\alpha$, $\alpha=1,2$, is a tangent one.

\noindent
To examine axially symmetric solutions to (1.1) we introduce the cylindrical 
coordinates $r$, $\varphi$, $z$ by the relations $x_1=r\cos\varphi$, 
$x_2=r\sin\varphi$, $x_3=z$.

\noindent
Moreover, we introduce the vectors
$$
\bar e_r=(\cos\varphi,\sin\varphi,0),\quad
\bar e_\varphi=(-\sin\varphi,\cos\varphi,0),\quad \bar e_z=(0,0,1)
$$
connected with the cylindrical coordinates.

\noindent
Then, the cylindrical coordinates of $v$ are defined by the relations
$$
v_r=v\cdot\bar e_r,\quad v_\varphi=v\cdot\bar e_\varphi,\quad
v_z=v\cdot\bar e_z,
$$
where the dot denotes the scalar product in $\R^3$.
Finally, by $u=rv_\varphi$ we denote a swirl.

To describe the domain $\Omega$ in greater details we introduce the notation
$$\eqal{
&\Omega=\{x\in\R^3:\ r<R,\ |z|<a\},\cr
&S_1=\{x\in\R^3:\ r=R,\ |z|<a\},\cr
&S_2=\{x\in\R^3,\ r<R,\ z\in\{-a,a\}\},\cr}
$$
where $R$ and $a$ are given positive numbers.

\Def{1.1.} 
By the axially symmetric solution we mean such solutions to problem (1.1) that
$$
v_{r,\varphi}=v_{\varphi,\varphi}=v_{z,\varphi}=p_{,\varphi}=0.
\leqno(1.2)
$$

In the cylindrical coordinates equations (1.1) for the axially symmetric 
solutions can be expressed in the form (see [LL, Ko])
$$
v_{r,t}+v\cdot\nabla v_r-{v_\varphi^2\over r}-\nu\Delta v_r+\nu
{v_r\over r^2}=-p_{,r},
\leqno(1.3)
$$
$$
v_{\varphi,t}+v\cdot\nabla v_\varphi+{v_r\over r}v_\varphi-\nu\Delta v_\varphi
+\nu{v_\varphi\over r^2}=0,
\leqno(1.4)
$$
$$
v_{z,t}+v\cdot\nabla v_z-\nu\Delta v_z=-p_{,z},
\leqno(1.5)
$$
$$
v_{r,r}+v_{z,z}=-{v_r\over r},
\leqno(1.6)
$$
where $v\cdot\nabla=v_r\partial_r+v_z\partial_z$, 
$\Delta u={1\over r}(ru_{,r})_{,r}+u_{,zz}$.
\goodbreak

\noindent
Let us introduce the $\varphi$-component of vorticity by
$$
\chi=v_{r,z}-v_{z,r}.
\leqno(1.7)
$$
Then $\chi$ is a solution to the problem (see [Z1])
$$\eqal{
&\chi_t+v\cdot\nabla\chi-{v_r\over r}\chi-\nu\bigg[\bigg(r\bigg(
{\chi\over r}\bigg)_{,r}\bigg)_{,r}+\chi_{,zz}+
2\bigg({\chi\over r}\bigg)_{,r}\bigg]={2v_\varphi v_{\varphi,z}\over r}\cr
&\chi|_{t=0}=\chi_0\cr
&\chi|_{S_1}=0,\ \ \chi|_{S_2}\ -\ {\rm periodic\ boundary\ conditon.}\cr}
\leqno(1.8)
$$

The existence of global regular axially symmetric 
solutions with vanishing swirl was proved by Ladyzhenskaya and 
Ukhovski-Yudovich, independently, long time ago (see [L2, UY]). 
Special cases as axially symmetric domains without the axis of symmetry 
or the case with sufficiently small swirl near the axis of symmetry at time 
$t=0$ were treated in [Z2, Z3], respectively. Moreover, global existence of 
large regular solutions with sufficiently small $v_\varphi$, $v_{r,\varphi}$, 
$v_{\varphi,\varphi}$, $v_{z,\varphi}$ at time $t=0$ were proved in [Z1, Z8] 
in cylindrical and axially symmetric domains, respectively. We have to add 
that the slip boundary conditions were necessary.

\noindent
Although a global $L_\infty$-estimate for swirl in the case of axially 
symmetric solutions was proved in [CL] (but under the technical assumption 
that $u$ vanishes at infinity and on the axis of symmetry) it was still not 
clear how to inrease regularity of such weak solutions. Having the existence 
of weak solutions to problem (1.1) such that $v\in V_2^0(\Omega^T)$ (see the 
definition of this space in Section 2) we obtain that 
$v\cdot\nabla v\in L_{5/4}(\Omega^T)$. Then by [S, ZZ] we obtain that 
solutions to (1.1) are such that $v\in W_{5/4}^{2,1}(\Omega^T)$, 
$\nabla p\in L_{5/4}(\Omega^T)$. Then $v\in L_{5/2}(\Omega^T)$, 
$\nabla p\in L_{5/3}(\Omega^T)$, so $v\cdot\nabla v\in L_1(\Omega^T)$ and the 
technique from [S, ZZ] can not be again applied.

\noindent
Passing to the cylindrical coordinates in estimates for the axially 
symmetric solution we obtain norms of weighted spaces (jacobian of 
transformation is $r$) so all imbeddings and interpolations are the same as 
in 3d case (see [Z3]).

\noindent
Examining regularity of the axially symmetric solutions in a neighborhood 
$\Omega_*$ located at a positive distance from the axis of symmetry we showed 
in [Z2] that $v'=(v_r,v_z)\in V_2^1(\Omega_*^T)$ (see Definition 2.10 in 
Section 2). Then by imbedding we have
$$
\|v'\|_{L_{10}(\Omega_*^T)}+\|\nabla v'\|_{L_{10/3}(\Omega_*^T)}\le c
\|v'\|_{V_2^1(\Omega_*^T)}
\leqno(1.9)
$$
(see Lemma 2.11 in Section 2).
\goodbreak

Then from (1.9) we have that $v'\cdot\nabla v\in L_2(\Omega_*^T)$ because 
$v'\in L_{10}(\Omega_*^T)$, $\nabla v\in L_2(\Omega_*^T)$, so by [S, ZZ] 
we obtain the existence of regular solutions such that 
$v\in W_{5/3}^{2,1}(\Omega_*^T)$, $\nabla p\in L_{5/3}(\Omega_*^T)$.
Hence $\nabla v\in L_{5/2}(\Omega_*^T)$. Then 
$v'\cdot\nabla v\in L_2(\Omega_*^T)$, so [S, ZZ] implies that 
$v\in W_2^{2,1}(\Omega_*^T)$, $\nabla p\in L_2(\Omega_*^T)$.

In this paper we prove (1.9) in a neighborhood of the axis of symmetry (see 
Lemma 6.4). Combining the estimates near the axis and far from the axis we obtain 
(1.9) in the whole domain $\Omega$.

\noindent
Then in [Z4] we prove the existence of regular solutions to problem (1.1) such 
that $v\in W_2^{2,1}(\Omega^T)$, $\nabla p\in L_2(\Omega^T)$.

\noindent
To formulate the main results of this paper we need the following assumptions

\noindent
{\bf Assumptions}

\item{1.} Assume that there exist positive finite constants $d_1$ and 
$d_2$ such that
$$
\|v_0\|_{L_2(\Omega)}\le d_1,
$$
and
$$
\|u_0\|_{L_\infty(\Omega)}\equiv d_2,
$$
where $u=rv_\varphi$, $u_0=u|_{t=0}$.

\item{2.} Assume that $u_0\in C^\alpha(\Omega)$, $\alpha\in(0,1)$.

\item{3.} Assume that $\zeta=\zeta_1(r)$ is a cutoff smooth function such that
$\zeta_1(r)=1$ for $r\le r_0$ and $\zeta_1(r)=0$ for $2r_0\le r<R$.\\
Assume that $\tilde v_\varphi=v_\varphi\zeta_1$, $\tilde\chi=\chi\zeta_1^2$ and
$$
\bigg\|{\tilde v_\varphi^2(0)\over r}\bigg\|_{L_2(\Omega)}<\infty,\quad
\bigg\|{\tilde\chi(0)\over r}\bigg\|_{L_2(\Omega)}<\infty.
$$

\item{4.} Assume that $r_0$ is so small that
$$
\|u\|_{L_\infty(\Omega_{\zeta_1}^T)}\le\root{4}\of{{5\over4}}\nu,
\leqno(1.10)
$$
where $\Omega_{\zeta_1}=\Omega\cap\supp\zeta_1$ (see assumptions of Lemma 5.3).
\vskip6pt

Now we describe restriction (1.10). We explain that it is not a smallness 
condition but restriction on $r_0$ only. In view of Assumption 2 and 
Lemma 2.7 we have that $u\in C^{\alpha,\alpha/2}(\Omega^T)$, where $\alpha$ 
is not larger that $1/2$. But Assumption 1 and Lemma 2.3 imply
$$
\bigg\|{u\over r^2}\bigg\|_{L_2(\Omega^T)}\le cd_1.
$$
Hence $u$ vanishes on the axis of symmetry.

\noindent
Therefore, since $u\in C(0,T;C^\alpha(\Omega))$ it follows that there exists 
$r_0$ so small that (1.10) can be satisfied.

\proclaim Theorem A. (see Lemmas 5.3, 5.5) 
Let the Assumptions 1--4 hold. Then for solutions to problem (1.1) the 
following a priori estimate is satisfied
$$\eqal{
&\bigg\|{\tilde v_\varphi\over r}\bigg\|_{L_4(\Omega^t)}^4+
\bigg\|{\tilde\chi\over r}\bigg\|_{V_2^0(\Omega^t)}^2+
\bigg\|\nabla\bigg({\tilde v_r\over r}\bigg)_{,r}\bigg\|_{L_2(\Omega^t)}^2\cr
&\quad+\bigg\|{1\over r}
\bigg({\tilde v_r\over r}\bigg)_{,r}\bigg\|_{L_2(\Omega^t)}^2\le\varphi
\bigg({1\over r_0},d_1,d_2,\nu\bigg)\bigg[1\cr
&\quad+\bigg\|{\tilde v_\varphi^2(0)\over r}\bigg\|_{L_2(\Omega)}+
\bigg\|{\tilde\chi(0)\over r}\bigg\|_{L_2(\Omega)}^2\bigg]\equiv A_0^2,\quad 
t\le T,\cr}
\leqno(1.11)
$$
where $\varphi$ is an increasing positive function and $\tilde v_r=v_r\zeta$ 
and $\zeta$ is introduced in Assumption 3.

\proclaim Theorem B. (see Lemma 6.4) 
Let the assumptions of Theorem A be satisfied. Then (see Definition 2.9)
$$
\|\tilde v'\|_{V_2^1(\Omega^t)}\le\varphi(A_0),\quad t\le T,
\leqno(1.12)
$$
where $v'=(v_r,v_z)$ and $\varphi$ is also an increasing positive function.

To prove Theorem B we need Assumption 4.
The assumption holds in view of Lemma 3.3. However, Lemma 3.3 
is proved under the existence of local solution to (1.1) such that 
$v\in W_2^{2,1}(\Omega^{T_*})$, $\nabla p\in L_2(\Omega^{T_*})$ with $T_*$ 
sufficiently small, the norm $\|u\|_{C^{\alpha,\alpha/2}}$ does not depend on 
the local solution. This is an important factor letting (1.12) for any $T$. 
In [Z4] there is shown how the local solution can be extended step by  step 
without its norm increasing. This way guarantees the existence of global 
regular solutions.

There are many results concerning the sufficient conditions of regularity of 
axially symmetric solutions. In [CL] such condition is\\ 
$\intop_0^Tdt\left(\intop_{\R^3}|v|^\gamma drd\varphi 
dz\right)^{\alpha/\gamma}<\infty$, where 
$1/\alpha+1/\gamma\le1/2$, $2<\gamma<+\infty$, $2<\alpha\le+\infty$. 
In [SZ] the following condition is assumed
$$
\esssup_{-1\le t\le0}\intop_\Omega{1\over r}|v(x,t)|^2dx<\infty.
$$
Finally in [SS] there are considered the following conditions
$$
\esssup_{Q(z_0,R)}r|\bar v(x,t)|<\infty
$$
and
$$
\esssup_{Q(z_0,R)}\sqrt{t_0-t}|\bar v(x,t)|<\infty
$$
where $z_0=(x_0,t_0)$, $Q(z_0,T)=B(x_0,R)\times(t_0-R^2,t_0)$, 
$B(x_0,R)=\{x\in\R^3:\ |x-x_0|<R\}$, 
$\bar v(x,t)=v_r\bar e_r+v_z\bar e_z$.

Now we outline the proofs of Theorems A and B. We have to emphasize that 
estimates (1.11) and (1.12) are a priori type estimates.

First we describe a proof of estimate (1.12). Lemma 2.1 yields the energy type 
estimate for weak solutions to (1.1). Next Lemma 2.2 yields the existence of 
local solution to problem (1.1) such that $v\in W_2^{2,1}(\Omega^{T_*})$, 
$\nabla p\in L_2(\Omega^{T_*})$, where $T_*$ is described by the assumptions 
of the lemma. Since we examine problem (1.1) in a bounded cylinder the 
$L_\infty$ bound for $u$ formulated in Lemma 2.5 is proved in Lemma 2.1 
in [Z5].

\noindent
Applying the DeGiorgi method developed by Ladyzhenskaya, Solonnikov, 
Uratseva (see [LSU, Ch. 2]) and assuming that $v\in W_2^{2,1}(\Omega^T)$ 
we show that $u\in C^{\alpha,\alpha/2}(\bar\Omega^T)$, where $\alpha=1/2$ 
(see Lemmas 2.7, 3.1 and [Z5]).

Having the H\"older continuity of $u$ Lemma 2.3 yields that $u$ vanishes on 
the axis of symmetry. The two properties are crucial to satisfy (1.10) for 
sufficiently small $r_0$. Then in a series of lemmas (see below) estimate 
(1.11) is proved for $t\le T_*$. Then elliptic problem (4.1) implies (1.12). 
The estimate gives regularity near the axis of symmetry because of properties 
of the cut-off function $\zeta=\zeta_1(r)$. From the proof of Lemma 5.3 in 
[Z4] we have that
$$
\tilde{\tilde v}'\in V_2^1(\Omega^{T_*}),
\leqno(1.13)
$$
where $\tilde{\tilde v}'=v'\zeta_2(r)$ and $\zeta_2(r)=0$ for 
$r\le{r_0\over2}$ and $\zeta_2(r)=1$ for $r\ge{3\over2}r_0$.

\noindent
Assuming that $\{\zeta_1(r),\zeta_2(r)\}$ compose a partition of unity in 
$\Omega$ imbeddings (1.12) and (1.13) imply that $v'\in V_2^1(\Omega^{T_*})$. 
This property yields that $v'(T_*)\in H^1(\Omega)$. To apply Lemma 2.2 for 
interval $[T_*,2T_*]$ we need also that $v_\varphi(T_*)\in H^1(\Omega)$, 
which is proved in Lemma 4.1 in [Z4].

\noindent
Moreover, we show in [Z4] that there exists a constant $\alpha_0$ and 
a corresponding to it $T_*>0$ such that if $\|v(0)\|_{H^1(\Omega)}\le\alpha_0$ 
then $\|v(T_*)\|_{H^1(\Omega)}\le\alpha_0$. This property implies that the 
local solution can be extended step by step.

Now we present shortly the steps of the proof of estimate (1.11). First we 
obtain (see Lemma 4.2)
$$\eqal{
&\bigg\|\nabla\bigg({\tilde v_r\over r}\bigg)_{,r}\bigg\|_{L_2(\Omega)}^2+
6\bigg\|{1\over r}\bigg({\tilde v_r\over r}\bigg)_{,r}\bigg\|_{L_2(\Omega)}^2
\le\bigg\|\bigg({\tilde\chi\over r}\bigg)_{,r}\bigg\|_{L_2(\Omega)}^2\cr
&\quad+\varphi\ {\rm (norms\ of\ data)},\cr}
\leqno(1.14)
$$
where $\varphi$ is an increasing positive function.

\noindent
Next Lemma 5.1 yields
$$
\bigg\|{\tilde\chi\over r}\bigg\|_{V_2^0(\Omega^t)}^2\le{1\over\nu}
\intop_{\Omega^t}{\tilde v_\varphi^4\over r^4}dxdt+\varphi\ {\rm (norms\ 
of\ data).}
\leqno(1.15)
$$
Finally, Lemma 5.2 implies the inequality
$$\eqal{
&{1\over4}\intop_\Omega{\tilde v_\varphi^4\over r^2}dx+{3\over4}\nu
\intop_{\Omega^t}\bigg|\nabla{\tilde v_\varphi^2\over r}\bigg|^2dxdt+
{1\over4}\nu\intop_{\Omega^t}{\tilde v_\varphi^4\over r^4}dxdt\cr
&\le\intop_{\Omega^t}\bigg|{v_r\over r}\bigg|{\tilde v_\varphi^4\over r^4}
dxdt+\varphi\ {\rm (norms\ of\ data).}\cr}
\leqno(1.16)
$$
From (1.14)--(1.16) and under the assumption (see (5.10))
$$
\|u\|_{L_\infty(\Omega_{\zeta_1}^t)}\le\root{4}\of{5\over4}\nu
\leqno(1.17)
$$
we obtain by Lemmas 5.3, 5.5 and Conclusion 5.6 the estimate (1.11).

The restriction (1.17) is satisfied in view of the H\"older continuity of $u$ 
and the property that $u$ vanishes on the axis of symmetry. Since 
$u\in C^{1/2,1/4}$ we can calculate precisely the number $r_0$ appearing 
in the definition of function $\zeta_1(r)$. This implies that the r.h.s. of 
(1.11) depends on $r_0$.

\section{2. Auxiliary results}

By $c$ we denote a generic constant which changes its value from line to 
line. A constant $c_k$ with index $k$ is defined by the first formula, where 
it appears. By $\varphi$ we denote the generic functions which changes its 
form from formula to formula and is always a positive and an increasing 
function.

\noindent
By $c(\sigma)$ we denote a generic constant which increases with $\sigma$.

Let us introduce the space
$$
V_2^0(\Omega^T)=\{u:\|u\|_{V_2^0(\Omega^T)}^2=\|u\|_{L_\infty(0,T;L_2(\Omega))}
+\nu\|\nabla u\|_{L_2(\Omega^T)}^2<\infty\}.
$$

\proclaim Lemma 2.1. 
Assume that $v_0\in L_2(\Omega)$. Then there exists a weak solution to problem 
(1.1) such that $v\in V_2^0(\Omega^T)$ and 
$$
\|v\|_{V_2^0(\Omega^T)}\le c_0\|v_0\|_{L_2(\Omega)}\equiv d_1.
\leqno(2.1)
$$
From properties of $V_2^0(\Omega^T)$ we have (see [LSU, Ch. 2, Sect. 3]
$$
\|v\|_{L_q(0,T;L_p(\Omega))}\le c_2\|v\|_{V_2^0(\Omega^T)},
\leqno(2.2)
$$
where
$$
{3\over p}+{2\over q}\ge{3\over2}.
$$

\proclaim Lemma 2.2. (see Theorem 3.1 in [Z4]). 
Assume that $v_0\in H^1(\Omega)$. Then for $T$ so small that
$$
c_*T^{1/2}\|v(0)\|_{H^1(\Omega)}\le1
$$
there exists a local solution to problem (1.1) such that 
$v\in W_2^{2,1}(\Omega^T)$, $\nabla p\in L_2(\Omega^T)$ and there exists 
a constant $c_0$ independent of $v$ and $p$ such that
$$
\|v\|_{W_s^{2,1}(\Omega^T)}+\|\nabla p\|_{L_2(\Omega^T)}\le
c_0\|v(0)\|_{H^1(\Omega)}.
\leqno(2.3)
$$

\proclaim Lemma 2.3. 
Let $v_0\in L_2(\Omega)$. Then the weak axially symmetric solutions to problem 
(1.1) satisfy the estimate
$$
\|v\|_{V_2^0(\Omega^T)}+\bigg\|{v_r\over r}\bigg\|_{L_2(\Omega^T)}^2+
\bigg\|{v_\varphi\over r}\bigg\|_{L_2(\Omega^T)}^2\le c_0^2
\|v_0\|_{L_2(\Omega)}^2\equiv d_1^2
\leqno(2.4)
$$

\noindent
The result is proved in Lemma 2.4 in [Z4].

\noindent
Let us introduce the quantity (swirl)
$$
u=rv_\varphi.
\leqno(2.5)
$$
We see that $u$ is a solution to the problem
$$\eqal{
&u_t+v\cdot\nabla u-\nu\Delta u+2\nu{u_{,r}\over r}=0\quad &{\rm in}\ \ 
\Omega^T,\cr
&u|_{t=0}=u_0,\quad &{\rm in}\ \ \Omega,\cr
&u_{,r}={2\over r}u\quad &{\rm on}\ \ S_1^T,\cr
&{\rm periodic\ boundary\ conditions}\quad &{\rm on}\ \ S_2^T,\cr
&\divv v=0\quad &{\rm in}\ \ \Omega,\cr}
\leqno(2.6)
$$
where the boundary condition $(2.6)_3$ is derived in [Z1, Ch. 4, Lemma 2.1].

\proclaim Lemma 2.4. (see Lemma 2.1 in [Z5]) 
Let $u$ be a solution to (2.6). Let $u_0\in L_\infty(\Omega)$. Then there 
exists $c$ independent of $u$ such that
$$
\|u\|_{L_\infty(\Omega^T)}\le c\|u_0\|_{L_\infty(\Omega)}\equiv d_2.
\leqno(2.7)
$$

\noindent
From (2.4) and (2.5) we derive the estimate
$$
\bigg\|{u\over r^2}\bigg\|_{L_2(\Omega^T)}\le d_1.
\leqno(2.8)
$$

\proclaim Lemma 2.5. 
Let the assumptions of Lemmas 2.3 and 2.4 hold. Then
$$
\|v_\varphi\|_{L_4(\Omega^T)}^4\le\|rv_\varphi\|_{L_\infty(\Omega^T)}^2
\bigg\|{v_\varphi\over r}\bigg\|_{L_2(\Omega^T)}^2\le d_1^2d_2^2.
\leqno(2.9)
$$

\noindent
In view of (2.7) we can prove

\proclaim Lemma 2.6. 
Assume that $u_0\in C^\alpha(\Omega_\zeta)$, where 
$\Omega_\zeta=\Omega\cap\supp\zeta$, $\zeta=\zeta(r)$ is smooth function such 
that $\zeta(r)=1$ for $r\le r_0$ and $\zeta(r)=0$ for $r\ge2r_0$, $2r_0<R$. 
Then any solution to (2.6) satisfies
$$
u\in C^{\alpha,\alpha/2}(\Omega_\zeta^T)
\leqno(2.10)
$$
where $\alpha$ equals $1/2$ (see [Z5]).

\noindent
We recall from [Z1, Ch. 2] the Hardy inequality

\proclaim Lemma 2.7. 
Assume that $f\in L_{p,-\mu}(\Omega)$ and $f_{,r}\in L_{p,1-\mu}(\Omega)$, 
where $L_{p,\nu}(\Omega)=\{u:\ \intop_\Omega|u|^pr^{p\nu}dx<\infty\}$, 
$p\in(1,\infty)$, $\nu\in\R$. Then
$$
\|f\|_{L_{p,-\mu}(\Omega)}\le{1\over\left|\mu-{2\over p}\right|}
\|f_{,r}\|_{L_{p,1-\mu}(\Omega)},
\leqno(2.11)
$$
where $\mu\not={2\over p}$.

\Remark{2.8.} (see [Z10]) 
If $\mu>{2\over p}$ then $f\not\in L_{p,-\mu}(\Omega)$ because 
$f(0)=f|_{r=0}\not=0$. In this case inequality (2.11) must be replaced by
$$
\|f-f(0)\|_{L_{p,-\mu}(\Omega)}\le
{1\over\left|\mu-{2\over p}\right|}\|f_{,r}\|_{L_{p,1-\mu}(\Omega)}.
\leqno(2.12)
$$

\Def{2.9.} 
We need also the space
$$\eqal{
V_2^k(\Omega^T)&=\{u:\ \|u\|_{V_2^k(\Omega^T)}=
\|u\|_{L_\infty(0,T;H^k(\Omega))}^2\cr
&\quad+\nu\|\nabla u\|_{L_2(0,T;H^k(\Omega))}^2<\infty\},\quad 
k\in\N\cup\{0\}.\cr}
\leqno(2.13)
$$
For $k=0$ we have the space $V_2^0(\Omega^T)$.

\noindent
Moreover
$$
H^k(\Omega)=\bigg\{u:\ \|u\|_{H^k(\Omega)}=\bigg(\sum_{|\alpha|\le k}
\intop_\Omega|D_x^\alpha u|^2dx\bigg)^{1/2}<\infty\bigg\},
\leqno(2.14)
$$
$k\in\N\cup\{0\}$. In [Z4] we prove

\proclaim Lemma 2.10. 
Let $\zeta'=\zeta_2(r)$ be a smooth positive function such that $\zeta'(r)=0$ 
for $r\le r_0/2$ and $\zeta'(r)=1$ for $r\ge r_0$. Let 
$\Omega_{\zeta'}=\Omega\cap\supp\zeta'$ and 
$\Omega_{\bar\zeta'}=\Omega\cap\{x\in\supp\zeta':\ \zeta'(r)=1\}$. 
Assume that $v_0\in L_2(\Omega)$ with $c_0\|v_0\|_{L_2(\Omega)}=d_1$. 
Assume also that $v_\varphi(0)\in L_{49/18}(\Omega_{\zeta'})$, 
$\chi(0)\in L_2(\Omega_{\zeta'})$. Then
$$
\|v'\|_{V_2^1(\Omega_{\bar\zeta'})}\le c(1/r_0,d_1)[d_1+
\|v_\varphi(0)\|_{L_{49/18}(\Omega_{\zeta'})}+
\|\chi(0)\|_{L_2(\Omega_{\zeta'})}],
\leqno(2.15)
$$
where $v'=(v_r,v_z)$.

\noindent
Repeating the considerations from [BIN, Ch. 2, Sect. 2.16] we have

\proclaim Lemma 2.11. 
Let $F(x)=\intop_0^xf(y)dy$, $\beta>{1\over p}$, 
$\|x^{-\beta+1}f(x)\|_{L_p(\R_\varepsilon)}<\infty$, $p\in(1,\infty)$, 
$\R_\varepsilon=\{x\in\R:\ 0<\varepsilon<x\}$. Let $f(x)=0$ for 
$x\in(0,\varepsilon]$. Then the following inequality holds
$$
\bigg(\intop_{\R_\varepsilon}|x^{-\beta}\intop_0^xf(y)dy|^pdx\bigg)^{1/p}\le
{1\over\beta-{1\over p}}\bigg(\intop_{\R_\varepsilon}
|x^{-\beta+1}f(x)|^pdx\bigg)^{1/p}.
\leqno(2.16)
$$

\Proof 
Making the change of variable and applying the generalized Min\-kow\-ski 
inequality (see [BIN, Ch. 2, Sect. 2.10]) we obtain
$$\eqal{
&\bigg(\intop_{\R_\varepsilon}
\bigg|x^{-\beta}\intop_0^xf(y)dy\bigg|^pdx\bigg)^{1/p}=
\bigg(\intop_{\R_\varepsilon}\bigg|x^{-\beta+1}\intop_0^1f(xt)dt\bigg|^p
dx\bigg)^{1/p}\cr
&\le\intop_0^1\bigg(\intop_{\R_\varepsilon}|x^{-\beta+1}f(xt)|^pdx\bigg)^{1/p}
dt\cr
&=\intop_0^1\bigg(\intop_{\R_{\varepsilon t}}|y^{-\beta+1}f(y)|^p
dy\bigg)^{1/p}t^{\beta-1-1/p}dt\cr
&=\intop_0^1t^{\beta-1-1/p}\bigg(\intop_{\R_\varepsilon}
|y^{-\beta+1}f(y)|^pdy\bigg)^{1/p}dt\equiv I,\cr}
$$
where we used that $f(y)=0$ for $y\le\varepsilon$.
Integrating with respect to $t$ yields (2.16). This concludes the proof.

Taking $u(x)=\intop_0^xf(y)dy$ and raising to the power $p$ the both sides of 
(2.16) we obtain
$$
\intop_{\R_\varepsilon}|r^{-\beta}u(r)|^pdr\le{1\over(\beta-1/p)^p}
\intop_{\R_\varepsilon}|r^{-\beta+1}u_{,r}|^pdr,
$$
where $u_{,r}=0$ for $r\le\varepsilon$.
\goodbreak

\noindent
Continuing,
$$
\intop_{\R_\varepsilon}|r^{-\beta-1/p}u(r)|^prdr\le{1\over(\beta-1/p)^p}
\intop_{\R_\varepsilon}|r^{-\beta-1/p+1}u_{,r}|^prdr.
\leqno(2.17)
$$
Assuming that $u=u(r,z)$ and integrating (2.17) with respect to $z$ yields
$$
\intop_\R\intop_{\R_\varepsilon}|r^{-\beta-1/p}u(r,z)|^pdx\le
{1\over(\beta-1/p)^p}\intop_\R\intop_{\R_\varepsilon}
|r^{-\beta-1/p+1}u_{,r}(r,z)|^pdx,
$$
where $dx=rdrdz$. Introducing $\alpha=\beta+1/p$ we obtain the following 
Hardy inequality
$$
\bigg(\intop_{\R\times\R_\varepsilon}|r^{-\alpha}u(r,z)|^pdx\bigg)^{1/p}\le
{1\over\alpha-2/p}\bigg(\intop_{\R\times\R_\varepsilon}
|r^{-\alpha+1}u_{,r}(r,z)|^pdx\bigg)^{1/p},
\leqno(2.18)
$$
where $\alpha>{2\over p}$, $u_{,r}(r,z)=0$ for $r\le\varepsilon$.

Using the cylindrical coordinates to examine problem (1.1) we derive 
equations with coefficients which are singular on the axis of symmetry (see 
(3.5), (3.6), (4.1) ,(4.11), (5.2)). To cancel the singularities we consider 
problem (1.1) in $\Omega_\varepsilon=\{x\in\Omega\colon 0<\varepsilon<r<R\}$. 
Therefore, we have to add some boundary conditions on the surface 
$r=\varepsilon$. We assume
$$
v\cdot\bar n|_{r=\varepsilon}=0,\quad 
\bar n\cdot\T(v,p)\cdot\bar\tau_\alpha|_{r=\varepsilon}=0,\quad \alpha=1,2,
\leqno(2.19)
$$
where $\bar n$ is the unit outward normal vector to the surface 
$r=\varepsilon$ so it is directed to the axis of symmetry.

\noindent
In view of [Z1, Ch. 4, Lemma 2.4] conditions (2.19) imply
$$
v_r|_{r=\varepsilon}=0,\quad v_{z,r}|_{r=\varepsilon}=0,\quad 
\bigg(v_{\varphi,r}+{1\over r}v_\varphi\bigg)\bigg|_{r=\varepsilon}=0,
\leqno(2.20)
$$
where the last equation assumes the form
$$
u_{,r}|_{r=\varepsilon}=0.
\leqno(2.21)
$$
Moreover, (2.20) implies
$$
\chi|_{r=\varepsilon}=0.
\leqno(2.22)
$$
Considering problem (1.1) in $\Omega_\varepsilon$ we should work with 
approximate functions $v_\varepsilon$, $u_\varepsilon$, $\chi_\varepsilon$ 
but we omit the index $\varepsilon$ for simplicity.

\noindent
Repeating the proof of Lemma 2.3 in the case $\Omega_\varepsilon$ we have

\proclaim Lemma 2.12. 
Let $v_0\in L_2(\Omega)$. Then
$$
\|v\|_{V_2^0(\Omega_\varepsilon^T)}+
\bigg\|{v_r\over r}\bigg\|_{L_2(\Omega_\varepsilon^T)}+
\bigg\|{v_\varphi\over r}\bigg\|_{L_2(\Omega_\varepsilon^T)}\le d_1.
\leqno(2.23)
$$

\section{3. Existence in $B_2(Q_T,M,\gamma,r,\delta,\varkappa)$. Idea of the 
proof of Lemma 2.6}

The H\"older continuity of $u$ in a neighborhood of the axis of symmetry is 
crucial in the proof of main estimate (5.28) (see also restriction (5.10) 
in Lemma 5.3).

In this Section we follow the ideas, definitions and results from [LSU, 
Ch. 2, Sects. 6, 7] which need some changes described in [Z5]. 
Since the result is very important we describe very roughly how the continuity 
follows. In this Section we base on the existence of weak solutions to 
problem (1.1).

\Def{3.1.} 
We say that $u\in V_2^0(Q_T)$ such that $\esssup_{Q_T}|u|\le M$ belongs to 
$B_2(Q_T\cap Q(\varrho,\tau),M,\gamma,r,\delta,\varkappa)$ if the function 
$w(x,t)=\pm u(x,t)$ satisfies the inequalities
$$\eqal{
&\max_{t_0\le t\le t_0+\tau}
\|w^{(k)}(x,t)\|_{L_2(B_{\varrho-\sigma_1\varrho})}^2\le
\|w^{(k)}(x,t_0)\|_{L_2(B_\varrho)}^2\cr
&\quad+\gamma[(\sigma_1\varrho)^{-2}\|w^{(k)}\|_{L_2(Q(\varrho,\tau))}^2+
\mu^{{2\over r}(1+\varkappa)}(k,\varrho,\tau)]\cr}
\leqno(3.1)
$$
and
$$\eqal{
&\|w^{(k)}\|_{V_2^0(Q(\varrho-\sigma_1\varrho,\tau-\sigma_2\tau))}^2\cr
&\le\gamma\{[(\sigma_1\varrho)^{-2}+(\sigma_2\tau)^{-1}]
\|w^{(k)}\|_{L_2(Q(\varrho,\tau))}^2+\mu^{{2\over r}(1+\varkappa)}
(k,\varrho,\tau)],\cr}
\leqno(3.2)
$$
where
$$
w^{(k)}(x,t)=\max\{w(x,t)-k;0\},
$$
$$\eqal{
Q(\varrho,\tau)&=B_\varrho(x_0)\times(t_0,t_0+\tau)\cr
&=\{(x,t)\in Q_T:\ 
|x-x_0|<\varrho,\ t_0<t<t_0+\tau\},\cr}
$$
$Q_T=\Omega\times(0,T)$, $\varrho,\tau$ are arbitrary positive numbers, 
$\sigma_1,\sigma_2$ -- arbitrary numbers from $(0,1)$,
$$\eqal{
&\mu(k,\varrho,\tau)=\intop_{t_0}^{t_0+\tau}{\rm meas}^{r\over q}A_{k,\varrho}
(t)dt,\cr
&A_{k,\varrho}(t)=\{x\in B_\varrho(x_0):\ w(x,t)>k\}.\cr}
$$
The numbers $M$, $\gamma$, $\delta$, $\varkappa$ are arbitrary positive 
and numbers $r$, $q$ satisfy the relation
$$
{2\over r}+{3\over q}={3\over2}.
\leqno(3.3)
$$
Finally, $k$ is a positive number satisfying the condition
$$
\esssup_{Q(\varrho,\tau)}w(x,t)-k\le\delta.
\leqno(3.4)
$$

\proclaim Lemma 3.2. 
Assume that $u$ satisfies (2.6) and $|u|\le M$. Assume that $\divv v=0$, 
$v\in L_{r'}(0,T;L_{q'}(\Omega))$, 
${3\over q'}+{2\over r'}=1-{3\over2}\varkappa$, $\varkappa>0$, $u$ 
is axially symmetric and periodic with respect to $x_3$. \\
Then $u^{(k)}\in B_2(Q_T\cap Q(\varrho,\tau),M,\gamma,r,\delta,\varkappa)$, 
where $M,k,\delta$ satisfy (3.4) in the form
$$
M-k\le\delta\ \ {\rm and}\ \ \varkappa={1\over6},\ \ r=q={10\over3}.
$$

\Proof (see the proof of Lemma 3.2 in [Z5]).

\proclaim Lemma 3.3.
Let $u$ be a solution to (2.6). \\
Let $\divv v'=0$, $v'\in L_{r'}(0,T;L_{q'}(B_\varrho))$ 
${2\over r'}+{3\over q'}=1-{3\over2}\varkappa$. Then 
$u\in C^{\alpha,\alpha/2}(Q_T\cap Q(\varrho,\tau))$, where 
$\alpha={3\varkappa\over2}$ and 
$\|u\|_{C^{\alpha,\alpha/2}}\le c\max\{2M,16\varrho_0^{3\varkappa\over2}\}$ 
with $c$ equal to some number.

\noindent
For $v\in W_2^{2,1}(B_\varrho^T)$ we have that $r'=q'=10$ so 
$\varkappa={1\over3}$ and $\alpha={1\over2}$.

\noindent
The Lemma is proved in Theorem 5.4 in [Z5].

\section{4. A priori estimates for $v_r$ and $v_z$ in a neighborhood of 
the $x_3$-axis}

In this section we examine the elliptic problem for $v_r=v_r(r,z,t)$, 
$v_z=v_z(r,z,t)$,
$$\eqal{
&v_{r,z}-v_{z,r}=\chi\quad &{\rm in}\ \ \Omega,\cr
&v_{r,r}+v_{z,z}+{v_r\over r}=0\quad &{\rm in}\ \ \Omega,\cr
&v_r|_{S_1}=0,\ {\rm periodicity\ with\ respect\ to}\ z.\cr}
\leqno(4.1)
$$
By localizing the problem to a neighborhood of the axis of symmetry we will 
be able to consider the result in $\R^3\cap\{z:\ |z|<a\}$ with periodicity 
with respect to $z$.

\noindent
Expressing $(4.1)_2$ in the form
$$
(rv_r)_{,r}+(rv_z)_{,z}=0
\leqno(4.2)
$$
we see that there exists a stream function $\psi$ such that
$$
v_r={\psi_{,z}\over r},\quad v_z=-{\psi_{,r}\over r}.
\leqno(4.3)
$$
Using (4.3) in $(4.1)_1$ we obtain the equation
$$
\bigg({\psi_{,z}\over r}\bigg)_{,z}+\bigg({\psi_{,r}\over r}\bigg)_{,r}=\chi.
\leqno(4.4)
$$
To estimate solutions of (4.4) near the axis of symmetry we introduce 
a cut-off smooth function $\zeta_1=\zeta_1(r)$ such that $\zeta_1(r)=1$ for 
$r\le r_0$ and $\zeta_1(r)=0$ for $r\ge2r_0$, where $2r_0<R$.

\noindent
Therefore, we introduce the notation
$$
\tilde\psi=\psi\zeta_1^2,\quad \tilde\chi=\chi\zeta_1^2.
\leqno(4.5)
$$
In virtue of (4.5) we express (4.4) in the form
$$
\bigg({\tilde\psi_{,z}\over r}\bigg)_{,z}+
\bigg({\tilde\psi_{,r}\over r}\bigg)_{,r}=\tilde\chi+{\psi_{,r}\over r}
(\zeta_1^2)_{,r}+\bigg({\psi(\zeta_1^2)_{,r}\over r}\bigg)_{,r}
\equiv\tilde\chi_*.
\leqno(4.6)
$$
For (4.6) we have the following boundary conditions
$$\eqal{
&\tilde\psi\ {\rm vanishes\ with\ all\ derivatives\ for}\ r\ge2r_0,\cr
&\tilde\psi\ {\rm is\ periodic\ with\ respect\ to}\ z\in[-a,a],\ {\rm so}\cr
&\tilde\psi(r,(2k+1)a,t)=\tilde\psi(r,(2k-1)a,t),\ \ k\in\Z.\cr}
\leqno(4.7)
$$
From [Z1--Z3] we know that a natural energy type estimate for solutions to 
(1.8) is an estimate for $\left\|{\chi\over r}\right\|_{V_2^0(\Omega^T)}$, 
because in this case the nonlinearity of the second and the third terms on 
the l.h.s. of (1.8) is eliminated. Hence we will be looking for estimates 
for $\tilde\psi$ in the case where the r.h.s. of (4.6) divided by $r$ 
belongs to $V_2^0(\Omega^T)$. This needs a very fast vanishing of $\tilde\psi$ 
and $\tilde\chi_*$ in a neighborhood of the axis of symmetry. To guarantee 
such vanishing we prove the existence of solutions to problem (4.6), (4.7) 
with such property. For this purpose we introduce a cylinder of radius 
$r=\varepsilon<r_0$. Then $\tilde\chi_*|_{r\le\varepsilon}=\chi$ and we 
assume that $\chi|_{r=\varepsilon}=0$. We take the idea from [L2]. The 
assumption is motivated by estimations presented in Section 5, where we have 
to work with functions ${\chi\over r}$, ${v_\varphi\over r}$. Then considering 
the problems for $\chi$ and $v_\varphi$ we have to add additional artificial 
boundary conditions $\chi|_{r=\varepsilon}=0$ and 
$v_\varphi|_{r=\varepsilon}=0$ (see problems (5.2) and (5.7)). In this way 
we get approximate functions $\chi_\varepsilon$ and $v_{\varphi\varepsilon}$ 
but we shall omit the index $\varepsilon$ for simplicity. Instead of the 
above approach we can treat the estimates in Section 4, 5, 6 as a priori 
estimates performed on smooth functions vanishing sufficiently fast near the 
axis of symmetry. However, we shall follow the first approach. To simplify 
problem (4.6), (4.7) we introduce the quantities $\eta$ and $\vartheta$ by 
the relations
$$
\tilde\psi=\eta r^2,\quad \tilde\chi_*=\vartheta r,
\leqno(4.8)
$$
where the restriction $\chi|_{r\le\varepsilon}$ implies that 
$\vartheta|_{r\le\varepsilon}=0$.

\noindent
Using (4.8) in (4.6) yields
$$
\eta_{,zz}+\eta_{,rr}+{3\eta_{,r}\over r}=\vartheta
\leqno(4.9)
$$
By the definition of $\Delta$ in the cylindrical coordinates we obtain
$$
\Delta\eta+{2\eta_{,r}\over r}=\vartheta.
\leqno(4.10)
$$
However, (4.10) contains a coefficient singular on the axis of symmetry we 
assume to consider (4.10) in domain $\Omega$, because the proof of Lemma 4.1 
basis on weak formulation of (4.10) which does not contain any singular 
term on the axis of symmetry.

\noindent
Then we examine the following boundary value problem
$$\eqal{
&\Delta\eta+{2\eta_{,r}\over r}=\vartheta\ \ {\rm in}\ \ \Omega,\cr
&\eta|_{r=2r_0}=0,\ {\rm periodic\ boundary\ conditions}\cr
&{\rm for}\ z\in\{(2k-1)a,(2k+1)a\},\ \ k\in\Z,\cr}
\leqno(4.11)
$$
where $\vartheta|_{r\le\varepsilon}=0$.

We restrict our considerations to the case $k=0$.

\proclaim Lemma 4.1. 
Assume that $\vartheta,\vartheta_{,r}\in L_2(\Omega)$, 
$\vartheta_{,\varphi}=0$, $\vartheta|_{r=2r_0}=0$,\\ $\eta_{,r}|_{r=2r_0}=0$. 
Then there exists a solution to problem (4.10), (4.11) such that 
$\eta_{,\varphi}=0$, $\eta\in H^1(\Omega)$, $\eta_{,r}\in H_0^1(\Omega)$, 
$\nabla\eta_{,z}\in L_2(\Omega)$, $\eta_{,zr}\in H_0^1(\Omega)$, where 
$\|u\|_{H_0^1(\Omega)}=\|\nabla u\|_{L_2(\Omega)}+
\left\|{u\over r}\right\|_{L_2(\Omega)}$ and the estimate
$$\eqal{
&\intop_\Omega(|\eta|^2+|\nabla\eta|^2+|\nabla\eta_{,r}|^2+
|\nabla\eta_{,z}|^2+|\nabla\eta_{,zr}|^2)dx\cr
&\quad+\intop_\Omega\bigg({\eta_{,r}^2\over r^2}+
{\eta_{,zr}^2\over r^2}\bigg)dx+\intop_{-a}^a(\eta^2+
\eta_{,z}^2+\eta_{,z}^2)|_{r=0}dz\cr
&\le c\intop_{\Omega_\varepsilon}(\vartheta^2+\vartheta_{,r}^2)dx.\cr}
\leqno(4.12)
$$
holds.

\Proof 
Multiplying $(4.11)_1$ by $\eta$, integrating over $\Omega$ and 
using boundary conditions $(4.11)_{2,3}$ we obtain
$$
\intop_{\Omega}(\eta_{,r}^2+\eta_{,z}^2)dx-
2\intop_{\Omega}\eta_{,r}\eta drdz=\intop_{\Omega}\vartheta\eta dx.
\leqno(4.13)
$$
The second term on the l.h.s. equals
$$
-\intop_{\Omega}(\eta^2)_{,r}drdz=\intop_{-a}^a\eta^2|_{r=0}dz.
$$
Hence, by the Poincare inequality we obtain from (4.13) the estimate
$$
\|\eta\|_{H^1(\Omega)}+\bigg(\intop_{-a}^a\eta^2|_{r=0}
dz\bigg)^{1/2}\le c_1\|\vartheta\|_{L_2(\Omega)},
\leqno(4.14)
$$
where $c_1$ is an increasing function of $r_0$.

\noindent
We prove existence of solutions to problem (4.11) by the Galerkin 
method. Looking for weak solutions in the form
$$
\eta^{(m)}=\sum_{i=1}^ma_{im}(t)\varphi_i(x),
$$
where $\{\varphi_i\}_{i=1}^\infty$ is a basis in $H^1(\Omega)$ 
and $a_{im}$, $i=1,\dots,m$, satisfy
$$\eqal{
&\sum_{i=1}^m\intop_{\Omega}a_{im}(t)\nabla'\varphi_i\nabla'
\varphi_jdx-2\sum_{i=1}^m\intop_{\Omega}a_{im}(t)\varphi_{i,r}\varphi_jdrdz\cr
&=\intop_{\Omega}\vartheta\cdot\varphi_jdx,\quad j=1,\dots,m,\cr}
\leqno(4.15)
$$
where $\nabla'=(\partial_r,\partial_z)$.
In view of (4.14) operator (4.11) is invertible so there exist 
solutions to (4.15) (see [LU, Ch. 3, Sect. 16]) and then by the Galerkin 
method we have the existence of solutions to problem (4.11) in $H^1(\Omega)$.

We increase regularity of the weak solution by getting estimates for higher 
derivatives. These estimates can be derived precisely by applying 
differences and passing to the appropriate limits. For simplicity we shall 
only restrict our considerations to show appropriate a priori estimates 
for the weak solutions.

\noindent
Differentiating $(4.11)_1$ with respect $r$, multiplying the result by 
$\eta_{,r}$ and integrating over $\Omega$ yields
$$
\intop_{\Omega}\Delta\eta_{,r}\eta_{,r}dx-3
\intop_{\Omega}{\eta_{,r}^2\over r^2}dx+2
\intop_{\Omega}{\eta_{,rr}\eta_{,r}\over r}dx=
\intop_{\Omega}\vartheta_{,r}\eta_{,r}dx,
$$
where we used that $(\Delta\eta)_{,r}=\Delta\eta_{,r}-{\eta_{,r}\over r^2}$.

\noindent
Integrating by parts and using that $\eta_{,r}|_{r=2r_0}=0$ we obtain
$$
\intop_{\Omega}|\nabla\eta_{,r}|^2dx+3\intop_{\Omega}
{\eta_{,r}^2\over r^2}dx-\intop_{\Omega}(\eta_{,r}^2)_{,r}drdz=
\intop_{\Omega}\vartheta(\eta_{,rr}r+\eta_{,r})drdz.
\leqno(4.16)
$$
The last term on the l.h.s. of (4.16) equals 
$\intop_{-a}^a\eta_{,r}^2|_{r=0}dz$ and the r.h.s. we estimate by
$$
\intop_{\Omega}\bigg({\varepsilon_1\over2}\eta_{,rr}^2+
{\varepsilon_2\over2}{\eta_{,r}^2\over r^2}\bigg)dx+
\bigg({1\over2\varepsilon_1}+{1\over2\varepsilon_2}\bigg)
\intop_{\Omega}\vartheta^2dx.
$$
Setting $\varepsilon_1=1$, $\varepsilon_2=5$, then (4.16) takes the form
$$
\intop_{\Omega}|\nabla\eta_{,r}|^2dx+\intop_{\Omega}{\eta_{,r}^2\over r^2}dx
+\intop_{-a}^a\eta_{,r}^2|_{r=0}dz\le{6\over5}\intop_{\Omega}\vartheta^2dx.
\leqno(4.17)
$$
Differentiating (4.10) with respect to $z$, multiplying by $\eta_{,z}$ 
and integrating over $\Omega$ yields
$$
\intop_{\Omega}\Delta\eta_{,z}\eta_{,z}dx+2
\intop_{\Omega}\eta_{,zr}\eta_{,z}drdz=
\intop_{\Omega}\vartheta_{,z}\eta_{,z}dx.
$$
Integrating by parts and using the periodicity condition gives
$$
\intop_{\Omega}|\nabla\eta_{,z}|^2dx+2\intop_{-a}^a
\eta_{,z}^2|_{r=0}dx\le\intop_{\Omega}\vartheta^2dx.
\leqno(4.18)
$$
Differentiating (4.10) with respect to $z$ yields
$$
\Delta\eta_{,z}+{2\eta_{,zr}\over r}=\vartheta_{,z}.
\leqno(4.19)
$$
Differentiating (4.19) with respect to $r$ gives
$$
\Delta\eta_{,zr}-{3\eta_{,zr}\over r^2}+2{\eta_{,zrr}\over r}=\vartheta_{,zr}.
\leqno(4.20)
$$
Multiplying (4.20) by $\eta_{,zr}$, integrating over $\Omega$ 
and using that $\eta_{,zr}|_{r=2r_0}=0$ and $(4.11)_{2,3}$ we obtain
$$\eqal{
&\intop_{\Omega}|\nabla\eta_{,zr}|^2dx+3
\intop_{\Omega}{\eta_{,zr}^2\over r^2}dx-
2\intop_{\Omega}{\eta_{,zrr}\eta_{,zr}\over r}dx\cr
&=-\intop_{\Omega}\vartheta_{,zr}\eta_{,zr}dx=
\intop_{\Omega}\vartheta_{,r}\eta_{,zzr}dx,\cr}
\leqno(4.21)
$$
where in view of the periodicity condition the integration by parts is 
performed in the r.h.s. of (4.21).

The last term on the l.h.s. of (4.21) equals
$$
-2\intop_{\Omega}\eta_{,zrr}\eta_{,zr}drdz=
-\intop_{\Omega}(\eta_{,zr}^2)_{,r}drdz=\intop_{-a}^a\eta_{,zr}^2|_{r=0}dz.
$$
Summarizing, we obtain from (4.21) the estimate
$$
\intop_{\Omega}|\nabla\eta_{,zr}|^2dx+6\intop_{\Omega}
{\eta_{,zr}^2\over r^2}dx+2\intop_{-a}^a\eta_{,zr}^2|_{r=0}dz\le
\intop_{\Omega}\vartheta_{,r}^2dx.
\leqno(4.22)
$$
From (4.14), (4.17), (4.18) and (4.22) we obtain (4.12). This concludes the 
proof.

\proclaim Lemma 4.2. 
Let $\zeta=\zeta_1(r)$ be a smooth function such that $\zeta_1(r)=1$ for 
$r\le r_0$ and $\zeta_1(r)=0$ for $r\ge2r_0$. Let $\tilde v_r=v_r\zeta_1^2$, 
$\tilde\chi=\chi\zeta_1^2$. Let 
$\left({{\tilde\chi\over r}}\right)_{,r}\in L_2(\Omega_\varepsilon)$, 
$v_z,v_{z,r}\in L_2(\Omega_\varepsilon)$. 
Then for sufficiently smooth solutions to problem (1.1) the following 
inequality is valid
$$\eqal{
&\bigg\|\nabla\bigg({\tilde v_r\over r}\bigg)_{,r}
\bigg\|_{L_2(\Omega_\varepsilon)}^2+
6\bigg\|{1\over r}\bigg({\tilde v_r\over r}\bigg)_{,r}
\bigg\|_{L_2(\Omega_\varepsilon)}^2\cr
&\le\bigg\|\bigg({\tilde\chi\over r}\bigg)_{,r}
\bigg\|_{L_2(\Omega_\varepsilon)}^2+c(1/r_0)
(\|v_z\|_{L_2(\Omega_\varepsilon)}^2\cr
&\quad+\|v_{z,r}\|_{L_2(\Omega_\varepsilon)}^2).\cr}
\leqno(4.23)
$$

\noindent
The inequality suggests fast decreasing of $\tilde v_r$ and $\tilde\chi$ 
approaching the axis of symmetry. This will be implied by the mentioned 
approximation and appropriate passing to the limit $\varepsilon=0$.

\Proof 
Using the definitions of $\tilde\chi_*$ by (4.6) and $\vartheta$ by (4.8) 
we obtain
$$
\vartheta={\tilde\chi\over r}+{1\over r}\bigg({\psi_{,r}\over r}
(\zeta_1^2)_{,r}+\bigg({\psi(\zeta_1^2)_{,r}\over r}\bigg)_{,r}\bigg)
\equiv{\tilde\chi\over r}+\vartheta_*.
$$
Hence we have
$$\eqal{
&\intop_\Omega\vartheta_{*,r}^2dx=\intop_\Omega\bigg[{1\over r}
\bigg({\psi_{,r}\over r}\zeta_{1,r}^2+\bigg({\psi\over r}
\zeta_{1,r}^2\bigg)_{,r}\bigg)\bigg]_{,r}^2dx\cr
&\le c\intop_\Omega\bigg[{1\over r^2}\bigg(2{\psi_{,r}\over r}\zeta_{1,r}^2+
\psi\bigg({\zeta_{1,r}^2\over r}\bigg)_{,r}\bigg)\bigg]^2dx\cr
&\quad+c\intop_\Omega
{1\over r^2}\bigg[\bigg({\psi_{,r}\over r}\bigg)_{,r}\zeta_{1,r}^2+
{\psi_{,r}\over r}\zeta_{1,rr}^2\cr
&\quad+\bigg({\psi\over r}\bigg)_{,rr}\zeta_{1,r}^2+2
\bigg({\psi\over r}\bigg)_{,r}\zeta_{1,rr}^2+
{\psi\over r}\zeta_{1,rrr}^2\bigg]^2dx\cr
&\le c\bigg({1\over r_0}\bigg)\intop_{\Omega\cap\supp\zeta_{1,r}}\bigg(v_z^2+
\bigg({\psi\over r^2}\bigg)^2+v_{z,r}^2\bigg)dx\cr
&\le c\bigg({1\over r_0}\bigg)\intop_{\Omega_\varepsilon}(v_z^2+v_{z,r}^2)dx,
\cr}
\leqno(4.24)
$$
where we used that $v_z={\psi_{,r}\over r}$ and 
$\supp\zeta_{1,r}=\{r:\ r_0\le r\le2r_0\}$.

\noindent
We employed also the Hardy inequality 
$\intop_\Omega{\psi^2\over r^4}dx\le\intop_\Omega{\psi_{,r}^2\over r^2}dx=
\intop_\Omega v_z^2dx$, where in view of (4.8) we have that 
$\psi_{,r}|_{r=0}=0$ and $\psi|_{r=0}=0$.

\noindent
In view of (4.24) and the expression for $\vartheta$ we obtain from (4.22) 
the inequality
$$\eqal{
&\intop_{\Omega_\varepsilon}|\nabla\eta_{,zr}|^2dx+6
\intop_{\Omega_\varepsilon}{\eta_{,zr}^2\over r^2}dx\cr
&\le\intop_{\Omega_\varepsilon}\bigg({\tilde\chi\over r}\bigg)_{,r}^2dx+
c\bigg({1\over r_0}\bigg)\intop_{\Omega_\varepsilon}(v_z^2+v_{z,r}^2)dx.\cr}
\leqno(4.25)
$$
Since $\eta={\tilde\psi\over r^2}$ we have that 
$\eta_{,z}={\tilde\psi_{,z}\over r^2}={\tilde v_r\over r}$. Then the second 
term on the l.h.s. of (4.25) takes the form
$$
\intop_{\Omega_\varepsilon}{1\over r^2}
\bigg({\tilde v_r\over r}\bigg)_{,r}^2dx.
\leqno(4.26)
$$
The first term on the l.h.s. of (4.25) equals
$$
\intop_{\Omega_\varepsilon}\bigg|\nabla\bigg({\tilde v_r\over r}
\bigg)_{,r}\bigg|^2dx.
\leqno(4.27)
$$
The term can not be estimated from below by (4.26) because the Hardy 
inequality (2.13) does not hold in this case. Therefore, (4.25) takes the form
$$\eqal{
&\intop_{\Omega_\varepsilon}\bigg|\nabla\bigg({\tilde v_r\over r}\bigg)_{,r}
\bigg|^2dx+6\intop_{\Omega_\varepsilon}{1\over r^2}
\bigg({\tilde v_r\over r}\bigg)_{,r}^2dx\le
\intop_{\Omega_\varepsilon}\bigg({\tilde\chi\over r}\bigg)_{,r}^2dx\cr
&\quad+c(1/r_0)\intop_{\Omega_\varepsilon}(v_z^2+v_{z,r}^2)dx.\cr}
\leqno(4.28)
$$
From (4.28) we obtain (4.23). This concludes the proof.

\section{5. Estimate}

Let us consider problem (1.8). Let $\zeta=\zeta_1(r)$ be the cut-off
function introduced in Section 4. Let
$$
\tilde\chi=\chi\zeta^2
\leqno(5.1)
$$
Then we introduce an approximate solution $\tilde\chi$ to problem (1.8) as 
a solution to the problem
$$\eqal{
&\tilde\chi_t+v\cdot\nabla\tilde\chi-{v_r\over r}\tilde\chi-\nu
\bigg[\bigg(r\bigg({\tilde\chi\over r}\bigg)_{,r}\bigg)_{,r}+\tilde\chi_{,zz}
+2\bigg({\tilde\chi\over r}\bigg)_{,r}\bigg]\cr
&=v\cdot\nabla\zeta^2\chi-\nu(\chi\zeta_{,r}^2)_{,r}-\nu r
\bigg({\chi\over r}\bigg)_{,r}\zeta_{,r}^2\cr
&\quad-2\nu\bigg({\chi\over r}\zeta_{,r}^2\bigg)+2
{\tilde v_\varphi\tilde v_{\varphi,z}\over r}\quad {\rm in}\ \ 
\Omega_\varepsilon,\cr
&\tilde\chi|_{t=0}=\tilde\chi_0,\ \ \tilde\chi|_{r=2r_0}=0, 
\tilde\chi|_{r=\varepsilon}=0\quad {\rm and}\ \ 
\tilde\chi|_{z=-a}=\tilde\chi|_{z=a}\cr}
\leqno(5.2)
$$
is periodic with respect to $z$, where $\tilde v_\varphi=v_\varphi\zeta$.

\proclaim Lemma 5.1. 
Let $\Omega_\varepsilon=\{x\in\R^3:\ 0<\varepsilon<r<R,\ |z|<a\}$. Assume that 
there exists a weak solution described by Lemma 2.3. Let 
${\tilde v_\varphi\over r}\in L_4(\Omega_\varepsilon^T)$. Let
$\chi\in L_{20\over7}(\Omega\cap\supp\zeta_{,r}\times(0,T))$, where 
$\supp\zeta_{,r}=\{r:\ r_0\le r\le2r_0\}$.
Then for sufficiently smooth solutions to problem (1.1) we have
$$\eqal{
&\bigg\|{\tilde\chi\over r}\bigg\|_{V_2^0(\Omega_\varepsilon\times(0,t))}^2\le 
c(1/r_0)d_1^2+c(1/r_0)d_1
\|\chi\|_{L_{{20\over7}(\Omega\cap\supp\zeta_{,r}\times(0,t))}}^2\cr
&\quad+\bigg\|{\tilde\chi(0)\over r}\bigg\|_{L_2(\Omega_\varepsilon)}^2+
{1\over\nu}\intop_0^t\intop_{\Omega_\varepsilon}
{\tilde v_\varphi^4\over r^4}dxdt,\quad t\le T,\cr}
\leqno(5.3)
$$
where $d_1$ is introduced in (2.1) and $\varepsilon<r_0$.

\Proof 
Multiplying (5.2) by ${\tilde\chi\over r^2}$ and integrating over 
$\Omega_\varepsilon$ we obtain
$$\eqal{
&{1\over2}{d\over dt}\bigg\|{\tilde\chi\over r}
\bigg\|_{L_2(\Omega_\varepsilon)}^2+\nu\bigg\|\nabla
\bigg({\tilde\chi\over r}\bigg)\bigg\|_{L_2(\Omega_\varepsilon)}^2\cr
&=\intop_{\Omega_\varepsilon}\bigg[v\cdot\nabla\zeta^2\chi
{\tilde\chi\over r^2}-\nu(\chi\zeta_{,r}^2)_{,r}{\tilde\chi\over r^2}-\nu r
\bigg({\chi\over r}\bigg)_{,r}\zeta_{,r}^2{\tilde\chi\over r^2}\cr
&\quad-2\nu\bigg({\chi\over r}\zeta_{,r}^2\bigg){\tilde\chi\over r^2}\bigg]dx+
2\intop_{\Omega_\varepsilon}{\tilde v_\varphi\tilde v_{\varphi,z}\over r}
{\tilde\chi\over r^2}dx.\cr}
\leqno(5.4)
$$
Now we estimate the particular terms in the r.h.s. of (5.4). The last term 
equals
$$\eqal{
&\intop_{\Omega_\varepsilon}
{(\tilde v_\varphi^2)_{,z}\over r^2}{\tilde\chi\over r}dx=
-\intop_{\Omega_\varepsilon}{\tilde v_\varphi^2\over r^2}
\bigg({\tilde\chi\over r}\bigg)_{,z}dx\le{\varepsilon_1\over2}
\intop_{\Omega_\varepsilon}\bigg({\tilde\chi\over r}\bigg)_{,z}^2dx\cr
&\quad+{1\over2\varepsilon_1}\intop_{\Omega_\varepsilon}
{\tilde v_\varphi^4\over r^4}dx.\cr}
$$
To estimate the terms in the first integral on the r.h.s. of (5.4) we use 
properties of the cut-off function $\zeta=\zeta(r)$.

\noindent
The first expression under the square bracked is bounded by
$$
c(1/r_0)\intop_{\Omega\cap\supp\zeta_{,r}}|v_r|\chi^2dx
$$
The fourth term under the square bracket is estimated by
$$
c(1/r_0)\intop_{\Omega\cap\supp\zeta_{,r}}\chi^2dx.
$$
The second term under the square bracket equals to
$$
-\nu\intop_\Omega\chi_{,r}\zeta_{,r}^2{\tilde\chi\over r^2}dx-\nu\intop_\Omega
\chi\zeta_{,rr}^2{\tilde\chi\over r^2}dx\equiv I_1.
$$
Integrating by parts in the first term in $I_1$ it takes the form
$$
-{\nu\over2}\intop_\Omega(\chi^2)_{,r}{(\zeta^2)_{,r}\zeta^2\over r^2}dx=
{\nu\over2}\intop_\Omega\chi^2\bigg({\zeta_{,r}^2\zeta^2\over r^2}\bigg)_{,r}
dx.
$$
Hence
$$
|I_1|\le c(1/r_0)\intop_{\Omega\cap\supp\zeta_{,r}}\chi^2dx.
$$
Similarly, the third term under the square bracket yields
$$
-{\nu\over2}\intop_\Omega\bigg(\bigg({\chi\over r}\bigg)^2\bigg)_{,r}
\zeta^2\zeta_{,r}^2dx={\nu\over2}\intop_\Omega{\chi^2\over r^2}
(\zeta_{,r}^2\zeta^2)_{,r}dx\equiv I_2,
$$
so
$$
|I_2|\le c(1/r_0)\intop_{\Omega\cap\supp\zeta_{,r}}\chi^2dx.
$$
Using the above estimates in (5.4) and assuming that $\varepsilon_1=\nu$ 
implies
$$\eqal{
&{1\over2}{d\over dt}\bigg\|{\tilde\chi\over r}
\bigg\|_{L_2(\Omega_\varepsilon)}^2+
{\nu\over2}\bigg\|\nabla{\tilde\chi\over r}
\bigg\|_{L_2(\Omega_\varepsilon)}^2\le
{1\over2\nu}\intop_{\Omega_\varepsilon}{\tilde v_\varphi^4\over r^4}dx\cr
&\quad+c(1/r_0)\intop_{\Omega\cap\supp\zeta_{,r}}\chi^2dx+c(1/r_0)
\intop_{\Omega\cap\supp\zeta_{,r}}|v_r|\chi^2dx.\cr}
\leqno(5.5)
$$
Integrating (5.5) with respect to time and using Lemma 2.3 we obtain
$$\eqal{
&\bigg\|{\tilde\chi\over r}\bigg\|_{V_2^0(\Omega_\varepsilon\times(0,t))}^2\le
{1\over\nu}\bigg\|{\tilde v_\varphi\over r}
\bigg\|_{L_4(\Omega_\varepsilon\times(0,t))}^4\cr
&\quad+c(1/r_0)d_1^2+c(1/r_0)\intop_0^t\intop_{\Omega\cap\supp\zeta_{,r}}
|v|\chi^2dxdt
+\bigg\|{\tilde\chi(0)\over r}\bigg\|_{L_2(\Omega_\varepsilon)}^2,\cr}
\leqno(5.6)
$$
where we used the definition of space $V_2^0(\Omega^T)$ from the beginning 
of Section~2.

\noindent
The third integral on the r.h.s. of (5.6) is estimated by
$$\eqal{
&c(1/r_0)\|v\|_{L_{{10\over3}(\Omega\cap\supp\zeta_{,r}\times(0,t))}}
\|\chi\|_{L_{{20\over7}(\Omega\cap\supp\zeta_{,r}\times(0,t))}}^2\cr
&\le c(1/r_0)d_1\|\chi\|_{L_{{20\over7}(\Omega\cap\supp\zeta_{,r}\times(0,t))}}^2,
\cr}
$$
where Lemma 2.3 was again used.
In view of above estimates we obtain from (5.6) inequality (5.3). This 
concludes the proof.
\kwadrat

Let us consider (1.4). Multiplying (1.4) by $\zeta=\zeta_1(r)$ and 
introducing the notation $\tilde v_\varphi=v_\varphi\zeta$ we assume an 
approximate $\tilde v_\varphi$ as a solution to the problem
$$\eqal{
&\tilde v_{\varphi,t}+v\cdot\nabla\tilde v_\varphi+{v_r\over r}\tilde v_\varphi
-\nu\Delta\tilde v_\varphi+\nu{\tilde v_\varphi\over r^2}=
v\cdot\nabla\zeta v_\varphi-2\nu\nabla v_\varphi\nabla\zeta\cr
&\quad-\nu v_\varphi\Delta\zeta\quad {\rm in}\ \ \Omega_\varepsilon,\cr
&\tilde v_\varphi|_{t=0}=\tilde v_\varphi(0),\quad 
\tilde v_\varphi|_{r=2r_0}=0,\quad 
\tilde v_\varphi|_{r=\varepsilon}=0,\quad
\tilde v_\varphi|_{z=-a}=\tilde v_\varphi|_{z=a},\cr}
\leqno(5.7)
$$
so $\tilde v_\varphi$ is periodic with respect to $z$

\proclaim Lemma 5.2. 
Assume that $\Omega_\varepsilon=\{x\in\R^3:\ 0<\varepsilon<r<R,\ |z|<a\}$, 
$\Omega_{\zeta_{,r}}=\Omega\cap\supp\zeta_{,r}$. Assume that $v$ is a weak 
solution to problem (1.1) satisfying (2.4), 
$rv_\varphi(0)\in L_\infty(\Omega)$, 
${\tilde v_\varphi^2(0)\over r}\in L_2(\Omega_\varepsilon)$. 
Assume that $v$ is sufficiently regular. 
Then the following a priori inequality holds
$$\eqal{
&{1\over4}\intop_{\Omega_\varepsilon}{\tilde v_\varphi^4\over r^2}dx+{3\over4}
\nu\intop_{\Omega_\varepsilon^t}
\bigg|\nabla{\tilde v_\varphi^2\over r}\bigg|^2dxdt+
{3\over4}\nu\intop_{\Omega_\varepsilon^t}{\tilde v_\varphi^4\over r^4}dxdt\cr
&\le{3\over2}\intop_0^t\intop_{\Omega_\varepsilon}\bigg|{v_r\over r}\bigg|
{\tilde v_\varphi^4\over r^2}dxdt+c\bigg({1\over r_0}\bigg)d_2^2(1+d_2)d_1^2+
{1\over4}\bigg\|{\tilde v_\varphi^2(0)\over r}
\bigg\|_{L_2(\Omega_\varepsilon)}^2,\cr}
\leqno(5.8)
$$
where $c(1/r_0)$ is an increasing function.

\Proof 
Multiplying (5.7) by ${\tilde v_\varphi|\tilde v_\varphi|^2\over r^2}$ and 
integrating over $\Omega_\varepsilon\times(0,t)$ we obtain
$$\eqal{
&{1\over4}\intop_{\Omega_\varepsilon}{\tilde v_\varphi^4\over r^2}dx+
{3\over4}\nu\intop_{\Omega_\varepsilon^t}
\bigg|\nabla{\tilde v_\varphi^2\over r}\bigg|^2dxdt+{3\over4}\nu
\intop_{\Omega_\varepsilon^t}{\tilde v_\varphi^4\over r^4}dxdt\cr
&\quad+\intop_{\Omega_\varepsilon^t}v\cdot\nabla\tilde v_\varphi
{\tilde v_\varphi|\tilde v_\varphi|^2\over r^2}dxdt\cr
&\le-\intop_0^t\intop_{\Omega_\varepsilon}{v_r\over r}
{\tilde v_\varphi^4\over r^2}dxdt-\intop_0^t\intop_{\Omega_\varepsilon}v\cdot
\nabla\zeta v_\varphi\tilde v_\varphi{\tilde v_\varphi^2\over r^2}dxdt\cr
&\quad-2\nu\intop_0^t\intop_{\Omega_\varepsilon}\nabla v_\varphi\nabla
\zeta\tilde v_\varphi{\tilde v_\varphi^2\over r^2}dxdt-\nu\intop_0^t
\intop_{\Omega_\varepsilon}v_\varphi
\Delta\zeta\tilde v_\varphi{\tilde v_\varphi^2\over r^2}dxdt\cr
&\quad+{1\over4}\bigg\|{\tilde v_\varphi^2(0)\over r}
\bigg\|_{L_2(\Omega_\varepsilon)}^2,\cr}
\leqno(5.9)
$$
where the expressions on the l.h.s. of the above inequality follow from the 
following calculations. The first, the fourth and the fifth terms from the 
l.h.s. of (5.7) imply
$$\eqal{
I&\equiv\intop_{\Omega_\varepsilon}\partial_t\tilde v_\varphi
{\tilde v_\varphi|\tilde v_\varphi|^2\over r^2}dx-\nu
\intop_{\Omega_\varepsilon}
\Delta\tilde v_\varphi{\tilde v_\varphi|\tilde v_\varphi|^2\over r^2}dx+
\nu\intop_{\Omega_\varepsilon}{\tilde v_\varphi^4\over r^4}dx\cr
&\equiv I_1+I_2+I_3,\cr}
$$
where
$$
I_1={1\over4}{d\over dt}\intop_{\Omega_\varepsilon}
{\tilde v_\varphi^4\over r^2}dx.
$$
Integrating by parts in $I_2$, using that $\tilde v_\varphi|_{r=2r_0}=0$, 
$\tilde v_\varphi|_{r=\varepsilon}=0$ and that $\tilde v_\varphi$ is periodic 
with respect to $z$ we have
$$\eqal{
I_2&=\nu\intop_{\Omega_\varepsilon}\nabla\tilde v_\varphi\cdot\nabla
\bigg({\tilde v_\varphi^3\over r^2}\bigg)dx=3\nu\intop_{\Omega_\varepsilon}
|\nabla\tilde v_\varphi|^2{\tilde v_\varphi^2\over r^2}dx-2\nu
\intop_{\Omega_\varepsilon}\tilde v_{\varphi,r}
{\tilde v_\varphi^3\over r^2}drdz\cr}
$$

$$\eqal{
&=3\nu\intop_{\Omega_\varepsilon}
\bigg|\nabla\tilde v_\varphi{\tilde v_\varphi\over r}\bigg|^2dx-{\nu\over2}
\intop_{\Omega_\varepsilon}{(\tilde v_\varphi^4)_{,r}\over r^2}drdz\cr
&={3\over4}\nu\intop_{\Omega_\varepsilon}
\bigg|{\nabla\tilde v_\varphi^2\over r}\bigg|^2dx-{\nu\over2}
\intop_{\Omega_\varepsilon}\bigg({\tilde v_\varphi^4\over r^2}\bigg)_{,r}drdz-
\nu\intop_{\Omega_\varepsilon}{v_\varphi^4\over r^4}dx\cr
&={3\over4}\nu\intop_{\Omega_\varepsilon}
\bigg|\nabla{\tilde v_\varphi^2\over r}+
{\tilde v_\varphi^2\over r^2}\nabla r\bigg|^2dx+{\nu\over2}\intop_{-a}^2
{\tilde v_\varphi^4\over r^2}\bigg|_{r=\varepsilon}dz\cr
&\quad-\nu\intop_{\Omega_\varepsilon}{v_\varphi^4\over r^4}dx=
{3\over4}\nu\intop_{\Omega_\varepsilon}
\bigg[\bigg|\nabla{\tilde v_\varphi^2\over r}\bigg|^2+
{\tilde v_\varphi^4\over r^4}+\partial_r\bigg({\tilde v_\varphi^2\over r}\bigg)
{\tilde v_\varphi^2\over r^2}\bigg]dx\cr
&\quad+{\nu\over2}\intop_{-a}^a
{\tilde v_\varphi^4\over r^2}\bigg|_{r=\varepsilon}dz-
\nu\intop_{\Omega_\varepsilon}{\tilde v_\varphi^4\over r^4}dx=
{3\over4}\nu\intop_{\Omega_\varepsilon}
\bigg|\nabla{\tilde v_\varphi^2\over r}\bigg|^2dx\cr
&\quad-{\nu\over4}\intop_{\Omega_\varepsilon}{\tilde v_\varphi^4\over r^4}dx+
{\nu\over2}\intop_{-a}^a{\tilde v_\varphi^4\over r^2}\bigg|_{r=\varepsilon}dz+
{3\over4}\nu\intop_{\Omega_\varepsilon}\partial_r
\bigg({\tilde v_\varphi^2\over r}\bigg){\tilde v_\varphi^2\over r}drdz\cr}
$$
The last term in $I_2$ equals
$$
{3\over8}\nu\intop_{\Omega_\varepsilon}\partial_r
\bigg({\tilde v_\varphi^4\over r^2}\bigg)drdz=-{3\over8}\nu\intop_{-a}^a
{\tilde v_\varphi^4\over r^2}\bigg|_{r=\varepsilon}dz
$$
Hence
$$
I_2={3\over4}\nu\intop_{\Omega_\varepsilon}
\bigg|\nabla{\tilde v_\varphi^2\over r}\bigg|^2dx-{\nu\over4}
\intop_{\Omega_\varepsilon}{\tilde v_\varphi^4\over r^4}dx+{\nu\over8}
\intop_{-a}^a{\tilde v_\varphi^4\over r^2}\bigg|_{r=\varepsilon}dz.
$$
Therefore $I$ takes the form
$$
I={1\over4}{d\over dt}\intop_{\Omega_\varepsilon}{\tilde v_\varphi^4\over r^2}
dx+{3\over4}\nu\intop_{\Omega_\varepsilon}
\bigg|\nabla{\tilde v_\varphi^2\over r}\bigg|^2dx+{3\over4}\nu
\intop_{\Omega_\varepsilon}{\tilde v_\varphi^4\over r^4}dx+{\nu\over8}
\intop_{-a}^a{\tilde v_\varphi^4\over r^2}\bigg|_{r=\varepsilon}dz.
$$
The last integral in $I$ vanishes. Integrating the result with respect to time 
implies the three integrals on the l.h.s. of (5.9).

\noindent
Finally, the last term on the l.h.s. of (5.9) equals
$$
{1\over2}\intop_{\Omega_\varepsilon^t}{v_r\over r}
{\tilde v_\varphi^4\over r^2}dxdt.
$$
The second integral on the r.h.s. of (5.9) we estimate by
$$\eqal{
&\intop_0^t\intop_{\Omega_{\zeta_{,r}}}|v|^2{|rv_\varphi|^3\over r^5}dxdt\le 
c(1/r_0)\|rv_\varphi\|_{L_\infty(\Omega\times(0,t))}^3\intop_0^t\intop_\Omega
|v|^2dxdt\cr
&\le c\bigg({1\over r_0}\bigg)d_2^3d_1^2.\cr}
$$
The third term on the r.h.s. of (5.9) takes the form
$$\eqal{
&-2\nu\intop_0^t\intop_\Omega
{\nabla v_\varphi v_\varphi|v_\varphi|^2\cdot\nabla\zeta\zeta|\zeta|^2
\over r^2}dxdt=-{\nu\over2}\intop_0^t\intop_\Omega
{\nabla(v_\varphi^4)\cdot\nabla\zeta|\zeta|^2\over r^2}dxdt\cr
&={\nu\over2}\intop_0^t\intop_\Omega v_\varphi^4\nabla\cdot
\bigg({\nabla\zeta\zeta|\zeta|^2\over r^2}\bigg)dxdt\equiv I_1.\cr}
$$
Hence,
$$
|I_1|\le c(1/r_0)\|rv_\varphi\|_{L_\infty(\Omega\times(0,t))}^2
\intop_{\Omega^t}v_\varphi^2dxdt\le c(1/r_0)d_2^2d_1^2.
$$

Similarly, the last but one term on the r.h.s. of (5.9) we estimate by 
$c\left({1\over r_0}\right)d_2^2d_1^2$.

\noindent
Using the above estimates in (5.9) implies (5.8). This concludes the proof.
\kwadrat

\proclaim Lemma 5.3. 
Assume that ${\tilde\chi(0)\over r}\in L_2(\Omega)$, 
${\tilde v_\varphi(0)\over r^{1/2}}\in L_4(\Omega)$, 
$\chi\in L_{{20\over7}(\Omega^T)}$. Let the assumptions of Lemmas 2.3, 2.4 
hold. Assume that
$$
\|rv_\varphi\|_{L_\infty(\Omega_\zeta^t)}\le\root{4}\of{5\over4}\nu,
\leqno(5.10)
$$
Then the following a priori inequality holds
$$\eqal{
&L^2(\Omega_\varepsilon^t)&\equiv\intop_{\Omega_\varepsilon^t}
{\tilde v_\varphi^4\over r^4}dxdt+
\bigg\|{\tilde\chi\over r}\bigg\|_{V_2^0(\Omega_\varepsilon^t)}^2+
\intop_{\Omega_\varepsilon^t}
\bigg|\nabla\bigg({\tilde v_r\over r}\bigg)_{,r}\bigg|^2dxdt\cr
&&\quad+\intop_{\Omega_\varepsilon^t}{1\over r^2}
\bigg({\tilde v_r\over r}\bigg)_{,r}^2dxdt
\le\varphi(d_1,d_2)
\|\chi\|_{L_{{20\over7}(\Omega\cap\supp\zeta_{,r}\times(0,t))}}^2\cr
&&\quad+\varphi(1/r_0,\nu,d_1,d_2)
\bigg(1+\bigg\|{\tilde v_\varphi^2(0)\over r}\bigg\|_{L_2(\Omega)}^2+
\bigg\|{\tilde\chi(0)\over r}\bigg\|_{L_2(\Omega)}^2\bigg),\cr
&t\le T,\cr}
\leqno(5.11)
$$
where $\tilde v_\varphi=v_\varphi\zeta$, $\tilde\chi=\chi\zeta^2$, 
$\tilde v_r=v_r\zeta^2$, $\zeta=\zeta_1(r)$, $\zeta_{1,r}\not=0$ for 
$r\in(r_0,2r_0)$, $2r_0<R$, 
$\Omega_\varepsilon=\{x\in\R^3:\ 0<\varepsilon<r<R,\ |z|<a\}$ and 
$\varphi$ is an increasing positive function. Moreover, 
$\Omega_\zeta=\Omega\cap\supp\zeta$.

\Proof 
From (4.28) we have
$$\eqal{
&\intop_{\Omega_\varepsilon}
\bigg|\nabla\bigg({\tilde v_r\over r}\bigg)_{,r}\bigg|^2dx+
6\intop_{\Omega_\varepsilon}{1\over r^2}
\bigg({\tilde v_r\over r}\bigg)_{,r}^2dx\le
\intop_{\Omega_\varepsilon}\bigg({\tilde\chi\over r}\bigg)_{,r}^2dx\cr
&\quad+c\intop_{\Omega_\zeta}(v_z^2+v_{z,r}^2)dx,\cr}
$$
where $\zeta=\zeta_1(r)$.

\noindent
Integrating the above inequality with respect to time, using estimate (5.3) 
and Lemma 2.3 in the r.h.s. we obtain
$$\eqal{
&\intop_{\Omega_\varepsilon^t}
\bigg|\nabla\bigg({\tilde v_r\over r}\bigg)_{,r}\bigg|^2dxdt+6
\intop_{\Omega_\varepsilon^t}{1\over r^2}\bigg({\tilde v_r\over r}\bigg)_{,r}^2
dxdt\le{1\over\nu^2}\intop_0^t\intop_{\Omega_\varepsilon}
{\tilde v_\varphi^4\over r^4}dxdt\cr
&\quad+{1\over\nu}c(1/r_0)d_1
\|\chi\|_{L_{20\over7}(\Omega\cap\supp\zeta_{,r}\times(0,t))}^2+
{1\over\nu}c(1/r_0)d_1^2\cr
&\quad+{1\over\nu}\bigg\|{\tilde\chi(0)\over r}\bigg\|_{L_2(\Omega)}^2.\cr}
\leqno(5.12)
$$
To estimate the first term on the r.h.s. of (5.12) we employ (5.8)
$$\eqal{
&\intop_{\Omega_\varepsilon\times(0,t)}{\tilde v_\varphi^4\over r^4}dxdt\le
{2\over\nu}\intop_0^t\intop_{\Omega_\varepsilon}\bigg|{v_r\over r}\bigg|
{\tilde v_\varphi^4\over r^2}dxdt\cr
&\quad+c\bigg({1\over r_0},\nu\bigg)d_2^2(1+d_2)d_1^2+{2\over\nu}
\bigg\|{\tilde v_\varphi^2(0)\over r}\bigg\|_{L_2(\Omega)}^2.\cr}
\leqno(5.13)
$$
Let us recall that $\tilde v_r=v_r\zeta^2$, $\tilde\chi=\chi\zeta^2$, 
$\tilde v_\varphi=v_\varphi\zeta$ and $\Omega_\zeta=\Omega\cap\supp\zeta$.

\noindent
The first term on the r.h.s. of (5.13) can be expressed in the form
$$\eqal{
&{2\over\nu}\intop_0^t\intop_{\Omega_\varepsilon}
\bigg|{v_r\zeta^2\over r}\bigg|
{v_\varphi^2\tilde v_\varphi^2\over r^2}dxdt
={2\over\nu}\intop_0^t\intop_{\Omega_\varepsilon}
\bigg|{\tilde v_r\over r^3}\bigg|
r^2v_\varphi^2{\tilde v_\varphi^2\over r^2}dxdt\cr
&\le{2\over\nu}\bigg[{\varepsilon\over2}\intop_0^t\intop_{\Omega_\varepsilon}
{\tilde v_\varphi^4\over r^4}dxdt+{1\over2\varepsilon}
\|rv_\varphi\|_{L_\infty(\Omega_\zeta^t)}^4\intop_0^t
\intop_{\Omega_\varepsilon}{\tilde v_r^2\over r^6}dxdt\bigg].\cr}
$$
Setting $\varepsilon=\nu/2$ and using the estimate in the r.h.s. of (5.13) 
yields
$$\eqal{
&\intop_{\Omega_\varepsilon^t}{\tilde v_\varphi^4\over r^4}dxdt\le{4\over\nu^2}
\|rv_\varphi\|_{L_\infty(\Omega_\zeta^t)}^4\intop_{\Omega^t}
{\tilde v_r^2\over r^6}dxdt\cr
&\quad+c\bigg({1\over r_0},\nu\bigg)d_2^2(1+d_2)d_1^2+{4\over\nu}
\bigg\|{\tilde v_\varphi^2(0)\over r}\bigg\|_{L_2(\Omega)}^2.\cr}
\leqno(5.14)
$$
Employing (5.14) in the r.h.s. of (5.12) we obtain
$$\eqal{
&6\intop_{\Omega_\varepsilon^t}{1\over r^2}
\bigg({\tilde v_r\over r}\bigg)_{,r}^2dxdt\le{4\over\nu^4}
\|rv_\varphi\|_{L_\infty(\Omega_\zeta^t)}^4\intop_{\Omega_\varepsilon^t}
{\tilde v_r^2\over r^6}dxdt\cr
&\quad+c\bigg({1\over r_0},\nu\bigg)d_2^2(1+d_2)d_1^2+
{1\over\nu}c(1/r_0)d_1^2\cr
&\quad+c(1/r_0)d_1\|\chi\|_{{20\over7}(\Omega_{\zeta_{,r}}\times(0,t))}^2\cr
&\quad+{4\over\nu}\bigg\|{\tilde v_\varphi^2(0)\over r}\bigg\|_{L_2(\Omega)}^2
\equiv{4\over\nu^4}\|rv_\varphi\|_{L_\infty(\Omega_\zeta^t)}^4
\intop_{\Omega_\varepsilon^t}
\bigg|{\tilde v_r\over r^3}\bigg|^2dxdt+A_1^2,\cr}
\leqno(5.15)
$$
where $c(1/r_0)$ is an increasing function.

\noindent
In view of the Hardy inequality (2.18) for $p=2$ inequality (5.15) takes 
the form
$$
6\intop_{\Omega_\varepsilon^t}{1\over r^2}
\bigg({\tilde v_r\over r}\bigg)_{,r}^2dxdt\le{4\over\nu^4}
\|rv_\varphi\|_{L_\infty(\Omega_\zeta^t)}^4\intop_{\Omega_\varepsilon}
{1\over r^2}\bigg({\tilde v_r\over r}\bigg)_{,r}^2dxdt+A_1^2.
\leqno(5.16)
$$
Assuming that
$$
6-{4\over\nu^4}\|rv_\varphi\|_{L_\infty(\Omega_\zeta^t)}^4\ge1
\leqno(5.17)
$$
we obtain from (5.16) the inequality
$$
\intop_{\Omega_\varepsilon^t}{1\over r^2}
\bigg({\tilde v_r\over r}\bigg)_{,r}^2dxdt\le A_1^2.
\leqno(5.18)
$$
Employing (5.18) in the r.h.s. of (5.14) yields
$$
\intop_{\Omega_\varepsilon^t}{\tilde v_\varphi^4\over r^4}dxdt\le
{2\over\nu^2}d_2^4A_1^2+cA_1^2
\leqno(5.19)
$$
Exploiting (5.19) in (5.3) and (5.12) implies
$$\eqal{
&\bigg\|{\tilde\chi\over r}\bigg\|_{V_2^0(\Omega_\varepsilon^t)}^2+
\intop_{\Omega_\varepsilon^t}
\bigg|\nabla\bigg({\tilde v_r\over r}\bigg)_{,r}\bigg|^2dxdt+
\intop_{\Omega_\varepsilon^t}
{1\over r^2}\bigg({\tilde v_r\over r}\bigg)_{,r}^2dxdt\cr
&\le c(1/r_0,\nu,d_1,d_2)A_1^2.\cr}
\leqno(5.20)
$$
This concludes the proof of the lemma.

\Remark{5.4.} 
Condition (5.10) looks like a smallness condition, however it can be justified 
by the following explanations.

\noindent
By Lemma 3.3 we have that $u=rv_\varphi\in C^{\alpha,\alpha/2}(\Omega^T)$, 
$\alpha=1/2$. Then (2.10) implies that $u|_{r=0}=0$.

Hence for a sufficiently small neighborhood of the $x_3$-axis, so for 
sufficiently small $r_0$, restriction (5.10) can be satisfied. Since 
$\alpha=1/2$ the number $r_0$ must be small but can be calculated precisely. 
Then a bound in the a priori estimate, which we are looking for, becomes an 
increasing function of $1/r_0$.

\noindent
The local existence is proved in interval $(0,T_*)$ and 
$v\in W_2^{2,1}(\Omega^{T_*})$. Next by Theorem B, Lemma 2.11 and Lemma 6.1 we 
are able to obtain a priori estimate for $v\in W_2^{2,1}(\Omega^T)$, 
$T>T_*$. Then $v(T_*)\in H^1(\Omega)$ and Lemma 2.2 yields existence in 
$W_2^{2,1}(\Omega\times(T_*,2T_*))$.

\noindent
Continuing, we have existence in $W_2^{2,1}(\Omega^T)$ for any $T>0$. This 
implies that in each time interval $(kT_*,(k+1)T_*)$, $k\in\N$, restriction 
(5.10) is satisfied.

To derive any estimate from (5.11) we have to estimate \\
$\|\chi\|_{{20\over7}(\Omega\times\supp\zeta_{,r}\times(0,T))}$. For this 
purpose we need

\proclaim Lemma 5.5. 
Let
$$
A_0=\varphi(1/r_0,\nu,d_1,d_2)\bigg(1+
\bigg\|{\tilde v_\varphi^2(0)\over r}\bigg\|_{L_2(\Omega)}^2+
\bigg\|{\tilde\chi(0)\over r}\bigg\|_{L_2(\Omega)}^2\bigg)
$$
be finite, where $\tilde v_\varphi=v_\varphi\zeta$, $\tilde\chi=\chi\zeta^2$ 
and the smooth cut-off function $\zeta$ is described in Lemma 5.3. Then there 
exists a positive constant $c$ such that
$$
\|\chi\|_{L_{{20\over7}(\Omega_\varepsilon\cap\supp\zeta_{,r}\times(0,T))}}
\le cA_0.
\leqno(5.21)
$$

\Proof 
From (5.11) we have
$$
\bigg\|{\tilde\chi\over r}\bigg\|_{L_{{10\over3}
(\Omega_\varepsilon\times(0,t))}}\le c_1
\bigg\|{\chi\over r}\bigg\|_{L_{{10\over3}(\Omega_\varepsilon\cap
\supp\zeta_{,r}\times(0,t))}}+A_0.
\leqno(5.22)
$$
Let us introduce the sets
$$
\Omega_\varepsilon^{(\lambda)}=\{(r,z):\ 0<\varepsilon<r\le 
r_0-\lambda,\ |z|<a\}
$$
and connect with them a set of cut-off smooth functions such that
$$
\zeta^{(\lambda)}=\cases{
1&for $(r,z)\in\Omega^{(\lambda)}_\varepsilon$\cr
0&for $(r,z)\in\Omega_\varepsilon\setminus\Omega_\varepsilon^{(\lambda/2)}$.
\cr}
$$
To simplify notation we denote $\eta={\tilde\chi\over r}$. Then (5.22) can be 
expressed in the form
$$
\|\eta\|_{L_{{10\over3}(\Omega_\varepsilon^{(\lambda)}\times(0,t))}}\le c_1
\|\eta\|_{L_{{10\over3}(\Omega_\varepsilon^{(\lambda/2)}\setminus
\Omega_\varepsilon^{(\lambda)}\times(0,t))}}+A_0.
\leqno(5.23)
$$
From (5.23) we have
$$
\intop_{\Omega_\varepsilon^{(\lambda)}\times(0,t)}|\eta|^{10\over3}dxdt\le c_2
c_1^{10\over3}
\intop_{\Omega_\varepsilon^{(\lambda/2)}\setminus\Omega_\varepsilon^{(\lambda)}\times(0,t)}
|\eta|^{10\over3}dxdt+c_2A_0^{10/3},
\leqno(5.24)
$$
where $c_2=2^{{10\over3}-1}$.

\noindent
Adding
$$
c_2c_1^{10\over3}\intop_{\Omega_\varepsilon^{(\lambda)}\times(0,t)}
|\eta|^{10\over3}dxdt
$$
to both sides of (5.24) we obtain
$$
\intop_{\Omega_\varepsilon^{(\lambda)}\times(0,t)}|\eta|^{10\over3}dxdt\le
{c_2c_1^{10\over3}\over1+c_2c_1^{10\over3}}
\intop_{\Omega_\varepsilon^{(\lambda/2)}\times(0,t)}|\eta|^{10\over3}dxdt+
{c_2\over1+c_2c_1^{10\over3}}A_0^{10\over3}.
\leqno(5.25)
$$
Introducing the notation
$$\eqal{
&f(\lambda)=\intop_{\Omega_\varepsilon^{(\lambda)}\times(0,t)}|\eta|^{10\over3}dxdt,
\quad
\mu={c_2c_1^{10\over3}\over1+c_2c_1^{10\over3}}<1,\cr
&K={c_2\over1+c_2c_1^{10/3}}A_0^{10\over3}\cr}
$$
we obtain from (5.25) the inequality
$$
f(\lambda)\le\mu f(\lambda/2)+K
\leqno(5.26)
$$
Hence we derive the estimate
$$
f(\lambda)\le\sum_{j=0}^\infty\mu^jK={1\over1-\mu}K
\leqno(5.27)
$$
Therefore, Lemma 5.5 is proved.

\Conclusion{5.6.}
In view of (5.21) inequality (5.11) implies the a priori estimate
$$
L^2(\Omega_\varepsilon^t)\le cA_0,\quad t\le T.
\leqno(5.28)
$$
where $L$ is defined by (5.11).

\Remark{5.7.}
All estimates in this section are made in 
$\Omega_\varepsilon=\{(r,z)\in\Omega:\ 0<\varepsilon<r<R,\ |z|<a\}$, 
$\varepsilon>0$. Then we introduced approximate functions $\chi_\varepsilon$, 
$v_{\varphi\varepsilon}$, $v_{r\varepsilon}$ which are such that
$$
v_{\varphi\varepsilon}|_{r=\varepsilon}=0, \quad
v_{r\varepsilon}|_{r=\varepsilon}=0,\quad
\chi_\varepsilon|_{r=\varepsilon}=0.
\leqno(5.29)
$$
They are defined by problems (5.2), (5.7) and (4.1), where the appropriate 
conditions from (5.29) are added. Then estimate (5.28) takes the form
$$
L_\varepsilon^2(\Omega_\varepsilon^t)\le cA_0,
\leqno(5.30)
$$
where $L_\varepsilon$ means that it depends on $\chi_\varepsilon$, 
$v_{\varphi\varepsilon}$, $v_{r\varepsilon}$. Then by appropriate passing 
with $\varepsilon$ to 0 we obtain (5.30) for $\varepsilon=0$. The passing 
will be performed in [Z4].

\section{6. Decay estimates}

In this section we derive estimates guaranteeing global existence.

The calculations in Lemma 6.1 are formal but they can be performed in the more 
rigorous way presented in Sections 4 and 5.

\proclaim Lemma 6.1. 
Let the assumptions of Lemmas 2.3 and 2.4 hold. Let\break
${\tilde v_\varphi^2(0)\over r}\in L_2(\Omega)$, 
${\tilde\chi(0)\over r}\in L_2(\Omega)$, where 
$\tilde v_\varphi=v_\varphi\zeta$, $\tilde\chi=\chi\zeta^2$ and 
$\zeta=\zeta(r)$ is introduced in Lemma 5.1. Let
$$
\|rv_\varphi\|_{L_\infty(\Omega_\zeta^T)}\le\root{4}\of{3\over4}\nu
$$
Then
$$\eqal{
&{1\over\nu^2}\bigg\|{\tilde v_\varphi^2(t)\over r}\bigg\|_{L_2(\Omega)}^2+
\bigg\|{\tilde\chi(t)\over r}\bigg\|_{L_2(\Omega)}^2\cr
&\le{1\over1-\varepsilon}[\varphi_1(1/r_0,d_1,d_2)+\varphi_2(1/r_0,d_1)\varepsilon^{-3}]d_1^2\cr
&\quad+{1\over1-\varepsilon}
\bigg({1\over\nu^2}\bigg\|{\tilde v_\varphi^2(0)\over r}\bigg\|_{L_2(\Omega)}^2
+\bigg\|{\tilde\chi(0)\over r}\bigg\|_{L_2(\Omega)}^2\bigg)
e^{-\nu_*t},\cr
&t\le T,\cr}
\leqno(6.1)
$$
where $\varphi_1$, $\varphi_2$ are increasing positive functions of their 
arguments, $\nu_*>0$ and $\varepsilon\in(0,1)$.

\Proof 
Let us introduce the notation
$$
\Omega_\zeta=\Omega\cap\supp\zeta,\quad 
\Omega_{\zeta_{,r}}=\Omega\cap\supp\zeta_{,r}.
$$
From (5.5) we have
$$\eqal{
&{d\over dt}\bigg\|{\tilde\chi\over r}\bigg\|_{L_2(\Omega)}^2+\nu
\bigg\|\nabla{\tilde\chi\over r}\bigg\|_{L_2(\Omega)}^2\cr
&\le{1\over\nu}\intop_\Omega{\tilde v_\varphi^4\over r^4}dx+c_1(1/r_0)
\intop_{\Omega_{\zeta_{r,}}}\chi^2dx+
c_2(1/r_0)\intop_{\Omega_{\zeta_{,r}}}|v|\chi^2dx,\cr}
\leqno(6.2)
$$
where $c_1(r_0)\sim r_0^{-4}$, $c_2(r_0)\sim r_0^{-3}$.

\noindent
Multiplying (5.7) ${\tilde v_\varphi|\tilde v_\varphi|^2\over r^2}$ and 
integrating over $\Omega$ we obtain
$$\eqal{
&{1\over4}{d\over dt}\bigg\|{\tilde v_\varphi^2\over r}\bigg\|_{L_2(\Omega)}^2
+{3\nu\over4}\bigg\|\nabla{\tilde v_\varphi^2\over r}\bigg\|_{L_2(\Omega)}^2
+{3\nu\over4}\intop_\Omega{\tilde v_\varphi^4\over r^4}dx\cr
&\le\|rv_\varphi\|_{L_\infty(\Omega_\zeta)}^2\intop_\Omega
{|\tilde v_r|\over r^3}{\tilde v_\varphi^2\over r^2}dx-\intop_\Omega
v\cdot\nabla\zeta v_\varphi\tilde v_\varphi{\tilde v_\varphi^2\over r^2}dx\cr
&\quad-2\nu\intop_\Omega\nabla v_\varphi\nabla\zeta\tilde v_\varphi
{\tilde v_\varphi^2\over r^2}dx-\nu\intop_\Omega v_\varphi\Delta\zeta
\tilde v_\varphi{\tilde v_\varphi^2\over r^2}dx.\cr}
\leqno(6.3)
$$
Now we examine the particular terms on the r.h.s. of (6.3). The second 
integral is bounded by
$$
{\varepsilon\over2}\intop_\Omega{\tilde v_\varphi^4\over r^4}dx+
{c_3(r_0)\over2\varepsilon}\intop_{\Omega_{\zeta_{,r}}}v^2v_\varphi^4dx,
$$
where $c_3(r_0)\sim r_0^{-2}$ and the second integral is estimated by
$$
{c_4(r_0)\over2\varepsilon}\|rv_\varphi\|_{L_\infty(\Omega)}^4
\intop_{\Omega_{\zeta_{r,}}}v^2dx,
$$
where $c_4(r_0)\sim r_0^{-6}$.

\noindent
Summarizing the second integral on the r.h.s. of (6.3) is bounded by
$$
{\varepsilon\over2}\intop_\Omega{\tilde v_\varphi^4\over r^4}dx+
{c_4(r_0)\over2\varepsilon}\|rv_\varphi\|_{L_\infty(\Omega)}^4
\intop_{\Omega_{\zeta_{,r}}}v^2dx.
\leqno(6.4)
$$
The third integral on the r.h.s. of (6.3) we express in the form
$$\eqal{
&-2\nu\intop_\Omega\nabla v_\varphi\nabla\zeta
{v_\varphi v_\varphi^2\zeta^3\over r^2}dx=-{\nu\over2}\intop_\Omega
\nabla v_\varphi^4\nabla\zeta{\zeta^3\over r^2}dx\cr
&={\nu\over2}\intop_\Omega v_\varphi^4\nabla\bigg(\nabla\zeta
{\zeta^3\over r^2}\bigg)dx={\nu\over2}\intop_\Omega v_\varphi^4\Delta\zeta
{\zeta^3\over r^2}dx+{3\over2}\nu\intop_\Omega v_\varphi^4|\nabla\zeta|^2
{\zeta^2\over r^2}dx\cr
&\quad-\nu\intop_\Omega v_\varphi^4\nabla\zeta\cdot\nabla r{\zeta^3\over r^3}
dx.\cr}
$$
Therefore the sum of the third and the fourth terms from the r.h.s. of (6.3) 
equals
$$
-{\nu\over2}\intop_\Omega v_\varphi^4\Delta\zeta{\zeta^3\over r^2}dx+
{3\over2}\nu\intop_\Omega v_\varphi^4|\nabla\zeta|^2{\zeta^2\over r^2}dx
-\nu\intop_\Omega v_\varphi^4\nabla\zeta\cdot\nabla r{\zeta^3\over r^3}dx
\equiv I_1.
$$
Hence
$$
|I_1|\le c_5(r_0)\intop_{\Omega_{\zeta_{,r}}}v_\varphi^4dx,
$$
where $c_5(r_0)\sim r_0^{-4}$.

\noindent
Assuming $\varepsilon={\nu\over2}$ in (6.4) and employing the above 
estimates in (6.3) we obtain
$$\eqal{
&{1\over4}{d\over dt}\bigg\|{\tilde v_\varphi^2\over r}\bigg\|_{L_2(\Omega)}^2
+{3\nu\over4}\bigg\|\nabla{\tilde v_\varphi^2\over r}\bigg\|_{L_2(\Omega)}^2
+{\nu\over2}\intop_\Omega{\tilde v_\varphi^4\over r^4}dx\cr
&\le\|rv_\varphi\|_{L_\infty(\Omega_\zeta)}^2\intop_\Omega
{|\tilde v_r|\over r^3}{\tilde v_\varphi^2\over r^2}dx+{c_4(r_0)\over\nu}
\|rv_\varphi\|_{L_\infty(\Omega)}^4\intop_{\Omega_{\zeta_{,r}}}v^2dx\cr
&\quad+c_5(r_0)\intop_{\Omega_{\zeta_{,r}}}v_\varphi^4dx.\cr}
\leqno(6.5)
$$
Applying the H\"older and the Young inequalities to the first term on the 
r.h.s. of (6.5) we have
$$\eqal{
&{1\over4}{d\over dt}\bigg\|{\tilde v_\varphi^2\over r}\bigg\|_{L_2(\Omega)}^2
+{3\nu\over4}\bigg\|\nabla{\tilde v_\varphi^2\over r}\bigg\|_{L_2(\Omega)}^2
+{\nu\over4}\intop_\Omega{\tilde v_\varphi^4\over r^4}dx\cr
&\le{1\over\nu}\|rv_\varphi\|_{L_\infty(\Omega_\zeta)}^4\intop_\Omega
{\tilde v_r^2\over r^6}dx+{c_4(r_0)\over\nu}\|rv_\varphi\|_{L_\infty(\Omega)}^4
\intop_{\Omega_{\zeta_{,r}}}v^2dx\cr
&\quad+c_5(r_0)\intop_{\Omega_{\zeta_{,r}}}v_\varphi^4dx.\cr}
\leqno(6.6)
$$
Using estimate (4.28) in (6.2) yields
$$\eqal{
&{d\over dt}\bigg\|{\tilde\chi\over r}\bigg\|_{L_2(\Omega)}^2+{\nu\over2}
\bigg\|\nabla{\tilde\chi\over r}\bigg\|_{L_2(\Omega)}^2+3\nu\intop_\Omega
{1\over r^2}\bigg({\tilde v_r\over r}\bigg)_{,r}^2dx\cr
&\quad+{\nu\over2}\intop_\Omega
\bigg|\nabla\bigg({\tilde v_r\over r}\bigg)_{,r}\bigg|^2dx
\le{1\over\nu}\intop_\Omega{\tilde v_\varphi^4\over r^4}dx+c_1(r_0)
\intop_{\Omega_{\zeta_{,r}}}\chi^2dx\cr
&\quad+c_2(r_0)\intop_{\Omega_{\zeta_{,r}}}
|v|\chi^2dx+c_6(r_0)\intop_{\Omega_{\zeta_{,r}}}(v_z^2+v_{z,r}^2)dx,\cr}
\leqno(6.7)
$$
where $c_6(r_0)\sim r_0^{-4}$.

\noindent
Multiplying (6.6) by ${4\over\nu^2}$, adding to (6.7) and applying the Hardy 
inequality to the first integral on the r.h.s. of (6.6) we obtain
$$\eqal{
&{d\over dt}\bigg({1\over\nu^2}
\bigg\|{\tilde v_\varphi^2\over r}\bigg\|_{L_2(\Omega)}^2+
\bigg\|{\tilde\chi\over r}\bigg\|_{L_2(\Omega)}^2\bigg)+{3\over\nu}
\bigg\|\nabla{\tilde v_\varphi^2\over r}\bigg\|_{L_2(\Omega)}^2
+{\nu\over2}\bigg\|\nabla{\tilde\chi\over r}\bigg\|_{L_2(\Omega)}^2\cr
&\quad+3\nu\intop_\Omega{1\over r^2}\bigg({v_r\over r}\bigg)_{,r}^2dx
\le{4\over\nu^3}\|rv_\varphi\|_{L_\infty(\Omega_\zeta)}^4\intop_\Omega
{1\over r^2}\bigg({v_r\over r}\bigg)_{,r}^2dx\cr
&\quad+{4c_4(r_0)\over\nu^3}
\|rv_\varphi\|_{L_\infty(\Omega)}^4\intop_{\Omega_{\zeta_{,r}}}v^2dx+
{4c_5(r_0)\over r^2}\intop_{\Omega_{\zeta_{,r}}}v_\varphi^4dx\cr
&\quad+c_1(r_0)\intop_{\Omega_{\zeta_{,r}}}\chi^2dx+c_2(r_0)
\intop_{\Omega_{\zeta_{,r}}}|v|\chi^2dx\cr
&\quad+c_6(r_0)\intop_{\Omega_{\zeta_{,r}}}(v_z^2+v_{z,r}^2)dx.\cr}
\leqno(6.8)
$$
Assuming
$$
3\nu\ge{4\over\nu^3}\|rv_\varphi\|_{L_\infty(\Omega_\zeta^T)}^4
\leqno(6.9)
$$
we obtain from (6.8) the inequality
$$\eqal{
&{d\over dt}\bigg({1\over\nu^2}
\bigg\|{\tilde v_\varphi^2\over r}\bigg\|_{L_2(\Omega)}^2+
\bigg\|{\tilde\chi\over r}\bigg\|_{L_2(\Omega)}^2\bigg)\cr
&\quad+\nu\bigg({3\over\nu^2}
\bigg\|\nabla{\tilde v_\varphi^2\over r}\bigg\|_{L_2(\Omega)}^2+{1\over2}
\bigg\|\nabla{\tilde\chi\over r}\bigg\|_{L_2(\Omega)}^2\bigg)
+{\nu\over2}\bigg\|\nabla{\tilde\chi\over r}\bigg\|_{L_2(\Omega)}^2\cr
&\le{4c_4(r_0)\over\nu^3}d_2^4\intop_{\Omega_{\zeta_{,r}}}v^2dx+
{4c_5(r_0)\over\nu^2}\intop_{\Omega_{\zeta_{,r}}}v_\varphi^4dx\cr
&\quad+c_1(r_0)\intop_{\Omega_{\zeta_{,r}}}\chi^2dx+
c_2(r_0)\intop_{\Omega_{\zeta_{,r}}}|v|\chi^2dx\cr
&\quad+c_6(r_0)\intop_{\Omega_{\zeta_{,r}}}(v_z^2+v_{z,r}^2)dx.\cr}
\leqno(6.10)
$$
Using the Poincare inequality in the second expression on the l.h.s. of 
(6.10) we can express (6.10) in the form
$$\eqal{
&{d\over dt}\bigg({1\over\nu^2}
\bigg\|{\tilde v_\varphi^2\over r}\bigg\|_{L_2(\Omega)}^2+
\bigg\|{\tilde\chi\over r}\bigg\|_{L_2(\Omega)}^2\bigg)\cr
&\quad+\nu_*\bigg({1\over\nu^2}\bigg\|{\tilde v_\varphi^2\over r}
\bigg\|_{L_2(\Omega)}^2+\bigg\|{\tilde\chi\over r}\bigg\|_{L_2(\Omega)}^2\bigg)
+{\nu\over2}\bigg\|\nabla{\tilde\chi\over r}\bigg\|_{L_2(\Omega)}^2\cr
&\le{4c_4(r_0)\over\nu^3}d_2^4\intop_{\Omega_{\zeta_{,r}}}v^2dx+
{4c_5(r_0)\over\nu^2}\intop_{\Omega_{\zeta_{,r}}}v_\varphi^4dx\cr
&\quad+c_1(r_0)\intop_{\Omega_{\zeta_{,r}}}\chi^2dx+c_2(r_0)
\intop_{\Omega_{\zeta_{,r}}}|v|\chi^2dx\cr
&\quad+c_6(r_0)\intop_{\Omega_{\zeta_{,r}}}(v_z^2+v_{z,r}^2)dx,\cr}
\leqno(6.11)
$$
where $\nu_*=\min\left\{{\nu\over2}{1\over c_p},{3\nu\over c_p}\right\}$ and 
$c_p$ is the constant from the Poincare inequality.

\noindent
Multiplying (6.11) by $e^{\nu_*t}$, integrating the result with respect to 
time we obtain
$$\eqal{
&\bigg({1\over\nu^2}
\bigg\|{\tilde v_\varphi^2(t)\over r}\bigg\|_{L_2(\Omega)}^2+
\bigg\|{\tilde\chi(t)\over r}\bigg\|_{L_2(\Omega)}^2\bigg)e^{\nu_*t}\cr
&\quad+{\nu\over2}\intop_0^t
\bigg\|\nabla{\tilde\chi(t')\over r}\bigg\|_{L_2(\Omega)}^2e^{\nu_*t'}dt'\cr
&\le c_7(r_0,d_1,d_2)d_1^2e^{\nu_*t}+c_2(r_0)\intop_0^t\intop_{\Omega_{\zeta_{,r}}}
|v|\,|\chi|^2dxe^{\nu_*t'}dt'\cr
&\quad+\bigg({1\over\nu^2}
\bigg\|{\tilde v_\varphi^2(0)\over r}\bigg\|_{L_2(\Omega)}^2+
\bigg\|{\tilde\chi(0)\over r}\bigg\|_{L_2(\Omega)}^2\bigg),\cr}
\leqno(6.12)
$$
where $c_7(r_0)\sim r_0^{-4}$.

\noindent
The second term on the r.h.s. is estimated by
$$\eqal{
&c_2(r_0)\intop_0^t\|v\|_{L_2(\Omega_{\zeta,r})}
\|\chi\|_{L_4(\Omega_{\zeta,r})}^2dt'\le\varepsilon\intop_0^t
\bigg\|\nabla{\chi\over r}\bigg\|_{L_2(\Omega_{\zeta,r})}^2dt'\cr
&\quad+c(1/r_0,d_1)\varepsilon^{-3}d_1^2,\cr}
\leqno(6.13)
$$
where $\varepsilon\in(0,1)$. Let
$$
X(t)={1\over\nu^2}\bigg\|{\tilde v_\varphi^2(t)\over r}\bigg\|_{L_2(\Omega)}^2
+\bigg\|{\tilde\chi(t)\over r}\bigg\|_{L_2(\Omega)}^2
\leqno(6.14)
$$
Employing (6.13) and (6.14) in (6.12) and applying the local interation 
argument (see [LSU, Ch. 4, Sect. 10]) we obtain
$$\eqal{
X(t)&\le{1\over1-\varepsilon}[c_7(r_0,d_1,d_2)d_1^2+c(1/r_0,d_1)
\varepsilon^{-3}d_1^2]\cr
&\quad+{1\over1-\varepsilon}X(0)e^{-\nu_*t},\cr}
\leqno(6.15)
$$
where $\varepsilon\in(0,1)$. This concludes the proof.

\noindent
Let us express inequality (6.15) in the form
$$
X(t)\le A+{1\over1-\varepsilon}X(0)e^{-\nu_*t},
\leqno(6.16)
$$
where $X(t)$ is defined by (6.14). Then we have

\proclaim Lemma 6.2. 
Let $T>0$ be so large that
$$
{1\over1-\varepsilon}e^{-\nu_*T}<1,\quad \varepsilon\in(0,1).
\leqno(6.17)
$$
Then for any $k\in\N$ the following estimate is valid
$$
X((k+1)T)\le{A\over1-1/(1-\varepsilon)e^{-\nu_*T}}+
{1\over1-\varepsilon}X(0)e^{-\nu_*kT}.
\leqno(6.18)
$$

\noindent
Proof follows easily by iteration.

\Remark{6.3.} 
In reality (6.1) is a global estimate, however, sometimes, it is convenient 
to examine problem (1.1) step by step in time. This may decrease a number 
of restrictions imposed on data (see [Z6]). Since (6.1) holds for any $T$ 
we see that (6.17) is not very restrictive. Moreover, $A$ is large because 
it is proportional to $r_0^{-4}\varepsilon^{-3}$, where $\varepsilon\in(0,1)$ 
and $r_0$ must be so small that the restriction holds
$$
\|u\|_{L_\infty(\Omega_\zeta^T)}\le\root{4}\of{3\over4}\nu.
$$
The above condition can be satisfied either by applying Lemma 3.3 and (2.8) 
with the H\"older exponent calculated in [Z5] or by applying the local in 
time regularity result (see Lemmas 2.2, 2.3). However, in the both cases 
$r_0$ must be small. We emphasize that smallness of $r_0$ implies that 
bounds either in (6.1) or in (6.18) are large.

\noindent
Now we shall show that the estimate
$$\eqal{
&\bigg\|{\tilde\chi\over r}\bigg\|_{V_2^0(\Omega^T)}\cr
&\qquad\le cA_0\equiv 
c\varphi(1/r_0,d_1,d_2)\bigg[1+
\bigg\|{\tilde v_\varphi^2(0)\over r}\bigg\|_{L_2(\Omega)}+
\bigg\|{\tilde\chi(0)\over r}\bigg\|_{L_2(\Omega)}\bigg]\cr}
\leqno(6.19)
$$
implies the following bound
$$
\|\tilde v'\|_{V_2^1(\Omega^T)}\le\varphi(A_0),
\leqno(6.20)
$$
where $v'=(v_r,v_z)$.

\noindent
For this purpose we come back to Section 4. Let us recall the notation
$$
\eta={\tilde\psi\over r^2},\quad \vartheta={\tilde\chi_*\over r},
$$
and $\eta$ is a solution to (4.10),
$$
\Delta\eta+{2\eta_{,r}\over r}=\vartheta,\quad
\Delta\eta={1\over r}(r\eta_{,r})_{,r}+\eta_{,zz}.
\leqno(6.21)
$$

\proclaim Lemma 6.4. 
Assume that (6.19) holds and $A_0$ is finite, where $\varphi$ is an increasing 
positive function, $d_1$ is defined by (2.4), $d_2$ by (2.7) and 
$\tilde v_\varphi(0)=v_\varphi(0)\zeta$, $\tilde\chi(0)=\chi(0)\zeta^2$. 
Moreover, $\zeta=\zeta(r)$ is a smooth cut off function such that 
$\zeta(r)=1$ for $r\le r_0$ and $\zeta(r)=0$ for $r\ge2r_0$. Finally, $r_0$ 
is so small that
$$
\|u\|_{L_\infty(\Omega_\zeta^T)}\le\root{4}\of{5\over4}\nu
$$
(see (5.10)), where $\Omega_\zeta=\Omega\cap\supp\zeta$. Then (6.20) holds.

\Proof 
Let us recall that $\eta$ and $\vartheta$ have compact supports with respect 
to $r$ and they are periodic with respect to $z$. From Section 4 formula 
(4.12) we have the estimate
$$\eqal{
&\intop_{\Omega_\varepsilon}(\eta^2+|\nabla\eta|^2+|\nabla\eta_{,r}|^2+
|\nabla\eta_{,z}|^2+|\nabla\eta_{,zr}|^2|)dx\cr
&\quad+\intop_{\Omega_\varepsilon}\bigg({\eta_{,r}^2\over r^2}+
{\eta_{,zr}^2\over r^2}\bigg)dx\le c\intop_{\Omega_\varepsilon}
(\vartheta^2+\vartheta_{,r}^2)dx\equiv I_1.\cr}
\leqno(6.22)
$$
From (6.22) we have
$$
\intop_{\Omega_\varepsilon}|\nabla\eta_{,zr}|^2dx=\intop_{\Omega_\varepsilon}
\bigg|\nabla\bigg({\tilde v_r\over r}\bigg)_{,r}\bigg|^2\le I_1
$$
so
$$
\intop_{\Omega_\varepsilon}
\bigg|\nabla\bigg({\tilde v_r\over r}\bigg)_{,r}\bigg|^2dx=
\intop_{\Omega_\varepsilon}\bigg[\bigg({\tilde v_r\over r}\bigg)_{,rr}^2+
\bigg({\tilde v_r\over r}\bigg)_{,rz}^2\bigg]dx\le I_1,
\leqno(6.23)
$$
where we used that $\eta={\tilde\psi\over r^2}$, $\tilde\psi=\psi\zeta_1^2$, 
$\tilde v_r=v_r\zeta_1^2$, $\tilde v_z=v_z\zeta_1^2$,
$$\eqal{
v_r&={\psi_{,z}\over r},\qquad v_z=-{\psi_{,r}\over r},\qquad
\eta_{,z}={\tilde\psi_{,z}\over r^2}={\tilde v_r\over r},\cr
\eta_{,r}&=\bigg({\tilde\psi\over r^2}\bigg)_{,r}=
{\psi_{,r}\over r^2}\zeta_1^2-{2\psi\over r^3}\zeta_1^2+
{\psi\over r^2}\zeta_{1,r}^2\cr
&=-{v_z\over r}\zeta_1^2-{2\psi\zeta_1^2\over r^3}+
{\psi\over r^2}\zeta_{1,r}^2=-{\tilde v_z\over r}-{2\tilde\psi\over r^3}+
{\psi\over r^2}\zeta_{1,r}^2.\cr}
\leqno(6.24)
$$
From (6.23) we have
$$
\intop_{\Omega_\varepsilon}\bigg({\tilde v_r\over r}\bigg)_{,rr}^2dx\le I_1,
$$
so
$$
\bigg({\tilde v_r\over r}\bigg)_{,rr}={\tilde v_{r,rr}\over r}-
2{\tilde v_{r,r}\over r^2}+2{\tilde v_r\over r^3}.
$$
Hence
$$
\intop_{\Omega_\varepsilon}{\tilde v_{r,rr}^2\over r^2}dx\le c
\bigg(\intop_{\Omega_\varepsilon}\bigg({\tilde v_r\over r}\bigg)_{,rr}^2dx+
\intop_{\Omega_\varepsilon}{\tilde v_{r,r}^2\over r^4}dx+
\intop_{\Omega_\varepsilon}{\tilde v_r^2\over r^6}dx\bigg)\equiv J_1.
$$
To estimate $J_1$ we calculate
$$\eqal{
&\intop_{\Omega_\varepsilon}{1\over r^2}
\bigg|{\tilde v_{r,r}\over r}\bigg|^2dx\cr
&\qquad=\intop_{\Omega_\varepsilon}{1\over r^2}
\bigg[\bigg({\tilde v_r\over r}\bigg)_{,r}-{\tilde v_r\over r^2}\bigg]^2dx\le 
c\bigg[\intop_{\Omega_\varepsilon}{1\over r^2}
\bigg({\tilde v_r\over r}\bigg)_{,r}^2dx+\intop_{\Omega_\varepsilon}
{\tilde v_r^2\over r^6}dx\bigg].\cr}
$$
Therefore, in view of (6.23), we obtain
$$
J_1\le cI_1+c\intop_{\Omega_\varepsilon}{1\over r^2}
\bigg({\tilde v_r\over r}\bigg)_{,r}^2dx+c\intop_{\Omega_\varepsilon}
{\tilde v_r^2\over r^6}dx\equiv J_2.
$$
To estimate $J_2$ we see that (6.22) implies
$$
\intop_{\Omega_\varepsilon}{\eta_{,zr}^2\over r^2}dx\le cI_1
$$
so
$$
\intop_{\Omega_\varepsilon}{1\over r^2}
\bigg({\tilde v_r\over r}\bigg)_{,r}^2dx\le cI_1.
$$
Moreover, by the Hardy inequality (2.18) for $p=2$, we have
$$
\intop_{\Omega_\varepsilon}{\tilde v_r^2\over r^6}dx=
\intop_{\Omega_\varepsilon}{1\over r^4}
\bigg|{\tilde v_r\over r}\bigg|^2dx\le\intop_{\Omega_\varepsilon}{1\over r^2}
\bigg({\tilde v_r\over r}\bigg)_{,r}^2dx\le cI_1.
$$
Using the above estimate in $J_2$ we obtain
$$
\intop_{\Omega_\varepsilon}{\tilde v_{r,rr}^2\over r^2}dx\le cI_1.
\leqno(6.25)
$$

Differentiating (6.21) twice with respect to $z$, multiplying the result by 
$\eta_{,zz}$ and integrating over $\Omega$ we obtain
$$
\intop_{\Omega}|\nabla\eta_{,zz}|^2dx-\intop_{\Omega}
(\eta_{,zz}^2)_{,r}drdz=\intop_{\Omega}\vartheta_{,z}\eta_{,zzz}dx,
$$
where the integration by parts was performed.

\noindent
Continuing, we have
$$
\intop_{\Omega}|\nabla\eta_{,zz}|^2dx+\intop_{-a}^a
\eta_{,zz}^2|_{r=0}dz\le c\intop_{\Omega}\vartheta_{,z}^2dx\equiv cI_2.
$$
Since
$$
\eta_{,zzz}=\bigg({\tilde v_r\over r}\bigg)_{,zz}={\tilde v_{r,zz}\over r}
$$
we obtain
$$
\intop_{\Omega_\varepsilon}{\tilde v_{r,zz}^2\over r^2}dx\le cI_2.
\leqno(6.26)
$$
We differentiate (6.21) twice with respect to $r$, multiply the result 
by $r^2\eta_{,rr}$ and integrate over $\Omega$. Then we obtain
$$
\intop_{\Omega}(\Delta\eta)_{,rr}r^2\eta_{,rr}dx+2\intop_{\Omega}
\bigg({\eta_{,r}\over r}\bigg)_{,rr}r^2\eta_{,rr}dx=\intop_{\Omega}
\vartheta_{,rr}r^2\eta_{,rr}dx.
\leqno(6.27)
$$
Using
$$
\Delta\eta=\eta_{,rr}+\eta_{,zz}+{\eta_{,r}\over r}
$$
we have
$$\eqal{
(\Delta\eta)_{,rr}&=\eta_{,rrrr}+\eta_{,rrzz}+{\eta_{,rrr}\over r}-
2{\eta_{,rr}\over r^2}+2{\eta_{,r}\over r^3}\cr
&=\Delta\eta_{,rr}-{2\eta_{,rr}\over r^2}+2{\eta_{,r}\over r^3}.\cr}
$$
Then (6.27) takes the form
$$\eqal{
&\intop_{\Omega}r^2\Delta\eta_{,rr}\eta_{,rr}dx-6\intop_{\Omega}
\eta_{,rr}^2dx+6\intop_{\Omega}{\eta_{,r}\eta_{,rr}\over r}dx\cr
&\quad+2\intop_{\Omega}\eta_{,rrr}\eta_{,rr}rdx=
\intop_{\Omega}\vartheta_{,rr}r^2\eta_{,rr}dx.\cr}
\leqno(6.28)
$$
Integrating by parts in the first term on the l.h.s. we obtain
$$\eqal{
&\intop_\Omega r^2\Delta\eta_{,rr}\eta_{,rr}dx=-\intop_\Omega\nabla
\eta_{,rr}\cdot\nabla(r^2\eta_{,rr})dx\cr
&=-\intop_\Omega|\nabla\eta_{,rr}|^2r^2dx-2\intop_\Omega\eta_{,rrr}
\eta_{,rr}rdx.\cr}
$$
Next,
$$\eqal{
&6\intop_\Omega{\eta_{,r}\eta_{,rr}\over r}dx=3\intop_\Omega(\eta_{,r}^2)_{,r}
drdz=-3\intop_{-a}^a\eta_{,r}^2|_{r=0}dz,\cr
&\intop_\Omega\vartheta_{,rr}r^2\eta_{,rr}dx=\intop_\Omega\vartheta_{,rr}
\eta_{,rr}r^3drdz=-\intop_\Omega\vartheta_{,r}(r^3\eta_{,rr})_{,r}drdz.\cr}
$$
In view of the above results equation (6.28) takes the form
$$\eqal{
&\intop_\Omega|\nabla\eta_{,rr}|^2r^2dx+6\intop_\Omega\eta_{,rr}^2dx+3
\intop_{-a}^a\eta_{,r}^2|_{r=0}dz\cr
&\quad+\intop_{-a}^a\eta_{,rr}^2r^2|_{r=0}dz\cr
&=\intop_\Omega\vartheta_{,r}(r^3\eta_{,rrr}+3r^2\eta_{,rr})drdz.\cr}
\leqno(6.29)
$$
Applying the H\"older and the Young inequalities in the r.h.s. of (6.29) and 
using that functions $\vartheta$ and $\eta$ have compact supports with respect 
to $r$, so $\vartheta$ and $\eta$ vanish for $r\ge2r_0$, we obtain
$$
\intop_\Omega|\nabla\eta_{,rr}|^2r^2dx+\intop_\Omega\eta_{,rr}^2dx\le
c(r_0)\intop_\Omega\vartheta_{,r}^2dx.
\leqno(6.30)
$$
Since $\eta={\tilde\psi\over r^2}$, (6.30) implies
$$
\intop_{\Omega_\varepsilon}\bigg|\nabla\bigg({\tilde\psi\over r^2}\bigg)_{,rr}
\bigg|^2r^2dx\le cI_1.
$$
Continuing, we see that
$$\eqal{
&\intop_{\Omega_\varepsilon}\bigg|\nabla\bigg({\tilde\psi\over r^2}\bigg)_{,rr}
\bigg|^2r^2dx=\intop_{\Omega_\varepsilon}\bigg|\nabla\bigg(
{\tilde\psi_{,r}\over r^2}-2{\tilde\psi\over r^3}\bigg)_{,r}\bigg|^2r^2dx\cr
&=\intop_{\Omega_\varepsilon}\bigg|\nabla\bigg({\tilde v_z\over r}+
{\psi\zeta_{1,r}^2\over r^2}-2{\tilde\psi\over r^3}\bigg)_{,r}\bigg|^2
r^2dx\le cI_1.\cr}
$$
Then
$$
\intop_{\Omega_\varepsilon}\bigg|\nabla\bigg({\tilde v_z\over r}\bigg)_{,r}
\bigg|^2r^2dx\le c\intop_{\Omega_\varepsilon}\bigg|\nabla
\bigg({\tilde\psi\over r^3}\bigg)_{,r}\bigg|^2r^2dx+
c\intop_{\Omega_\varepsilon}\bigg|\nabla{\psi\zeta_{1,r}^2\over r^2}\bigg|^2
r^2dx+cI_1.
\leqno(6.31)
$$
To examine the first integral on the r.h.s. of (6.31) we recall that 
$\tilde\psi=\eta r^2$. Then
$$\eqal{
&\intop_{\Omega_\varepsilon}\bigg|\nabla\bigg({\tilde\psi\over r^3}\bigg)_{,r}
\bigg|^2r^2dx=\intop_{\Omega_\varepsilon}\bigg|\nabla
\bigg({\eta\over r}\bigg)_{,r}\bigg|^2r^2dx=\intop_{\Omega_\varepsilon}
\bigg|\nabla\bigg({\eta_{,r}\over r}-{\eta\over r^2}\bigg)\bigg|^2r^2dx\cr
&\le c\intop_{\Omega_\varepsilon}\bigg(\bigg|\nabla{\eta_{,r}\over r}\bigg|^2+
\bigg|\nabla{\eta\over r^2}\bigg|^2\bigg)r^2dx\cr
&\le c\intop_{\Omega_\varepsilon}|\nabla\eta_{,r}|^2dx+
c\intop_{\Omega_\varepsilon}\bigg({\eta_{,r}^2\over r^2}+
{|\nabla\eta|^2\over r^2}+{\eta^2\over r^4}\bigg)dx\equiv K_1.\cr}
$$
In view of (6.22) we obtain
$$
K_1\le cI_1+c\intop_{\Omega_\varepsilon}\bigg({\eta_{,z}^2\over r^2}+
{\eta^2\over r^4}\bigg)dx.
$$
since the support of $\eta$ with respect to $r$ is compact, because $\eta$ 
vanishes for $r\ge2r_0$, we have
$$\eqal{
&\intop_{\Omega_\varepsilon}{\eta_{,z}^2\over r^2}dx\le c(r_0)
\intop_{\Omega_\varepsilon}{\eta_{,z}^2\over r^{2+2\delta}}dx\le c(r_0)
\intop_{\Omega_\varepsilon}{\eta_{,zr}^2\over r^{2\delta}}dx\cr
&\le c(r_0)\intop_{\Omega_\varepsilon}{\eta_{,zr}^2\over r^2}dx\le c(r_0)I_1,
\cr}
$$
where $\delta\in(0,1)$.

\noindent
By the Hardy inequality (2.18) for $p=2$ with $\eta_{,r}=0$ for 
$r\le\varepsilon$ we obtain
$$
\intop_{\Omega_\varepsilon}{\eta^2\over r^4}dx\le c\intop_{\Omega_\varepsilon}
{\eta_{,r}^2\over r^2}dx\le cI_1.
$$
Hence,
$$
K_1\le cI_1
$$
and (6.31) assumes the form
$$
\intop_{\Omega_\varepsilon}\bigg|\nabla\bigg({\tilde v_z\over r}\bigg)_{,r}
\bigg|^2r^2dx\le cI_1+c\intop_{\Omega_\varepsilon}\bigg|\nabla
{\psi\zeta_{1,r}^2\over r^2}\bigg|^2r^2dx.
\leqno(6.32)
$$
Since
$$
\nabla\bigg({\tilde v_z\over r}\bigg)_{,r}={\nabla\tilde v_{z,r}\over r}-
{\tilde v_{z,r}\over r^2}-{\nabla\tilde v_z\over r^2}+2{\tilde v_z\over r^3}
\quad {\rm in}\ \ \Omega_\varepsilon
$$
we obtain from (6.32) the inequality
$$\eqal{
&\intop_{\Omega_\varepsilon}|\nabla\tilde v_{z,r}|^2dx\le cI_1+
c\intop_{\Omega_\varepsilon}\bigg(\bigg|{\tilde v_{z,r}\over r^2}\bigg|^2+
\bigg|{\tilde v_{z,z}\over r^2}\bigg|^2+{\tilde v_z^2\over r^6}\bigg)r^2dx\cr
&\quad+c\intop_{\Omega_\varepsilon}\bigg|\nabla
{\psi\zeta_{1,r}^2\over r^2}\bigg|^2r^2dx\cr
&\le cI_1+c\intop_{\Omega_\varepsilon}\bigg({\tilde v_{z,r}^2\over r^2}+
{\tilde v_{z,z}^2\over r^2}+{\tilde v_z^2\over r^4}\bigg)dx+c
\intop_{\Omega_\varepsilon}\bigg|\nabla
\bigg({\psi\zeta_{1,r}^2\over r^2}\bigg)\bigg|^2r^2dx\cr}
\leqno(6.33)
$$
Now we shall estimate the middle integral on the r.h.s. of (6.33). Since 
$v_z=-{\psi_{,r}\over r}$ and $\tilde\psi=\eta r^2$ we have that
$$\eqal{
&\tilde v_{z,z}=-\eta_{,rz}r-2\eta_{,z}+{\psi_{,z}\over r}\zeta_{1,r}^2\cr
&\tilde v_z=-\eta_{,r}r-2\eta+{\psi\over r}\zeta_{1,r}^2,\quad
v_{z,r}=-\eta_{,rr}r-3\eta_{,r}+\bigg({\psi\over r}\zeta_{1,r}^2\bigg)_{,r},
\cr}
$$
Then, in view of (6.22), we have
$$\eqal{
\intop_{\Omega_\varepsilon}{\tilde v_{z,r}^2\over r^2}dx
&\le c\intop_{\Omega_\varepsilon}(\eta_{,rr}^2+{\eta_{,r}^2\over r^2}+
{1\over r^2}\bigg({\psi\over r}\zeta_{1,r}^2\bigg)_{,r}^2\bigg)dx\cr
&\le cI_1+c\bigg({1\over r_0}\bigg)\intop_{\Omega_\varepsilon}
\bigg({\psi\over r}\zeta_{1,r}^2\bigg)_{,r}^2dx,\cr
\intop_{\Omega_\varepsilon}{\tilde v_{z,z}^2\over r^2}dx
&\le c\intop_{\Omega_\varepsilon}\bigg(\eta_{,rz}^2+{\eta_{,z}^2\over r^2}+
v_r^2\zeta_{1,r}^2\bigg)dx\le cI_1+c(1/r_0)\intop_{\Omega_\varepsilon}
v_r^2dx,\cr
\intop_{\Omega_\varepsilon}{\tilde v_z^2\over r^4}dx
&\le c\intop_{\Omega_\varepsilon}\bigg({\eta_{,r}^2\over r^2}+
{\eta^2\over r^4}\bigg)dx+c(1/r_0)\intop_{\Omega_\varepsilon}\psi^2dx.\cr}
$$
In view of (6.22) and the Hardy inequality we have
$$
\intop_{\Omega_\varepsilon}\bigg({\eta_{,r}^2\over r^2}+
{\eta^2\over r^4}\bigg)dx\le cI_1.
$$
Since $v_r|_{r=0}=0$ we have that $\psi_{,z}|_{r=0}=0$, so $\psi|_{r=0}=0$ 
can be chosen. Then
$$
\psi=\intop_0^r\psi_{,r}dr=\intop_0^r{\psi_{,r}\over r}rdr
$$
and
$$
\intop_\Omega\psi^2dx\le c\intop_\Omega v_z^2dx.
\leqno(6.34)
$$
Moreover,
$$
\intop_{\Omega_\varepsilon}\bigg({\psi\over r}\zeta_{1,r}^2\bigg)_{,r}^2dx\le
c(1/r_0)\bigg(\intop_{\Omega_\varepsilon}{\psi_{,r}^2\over r^2}dx+
\intop_\Omega\psi^2dx\bigg)\le c(1/r_0)\intop_\Omega v_z^2dx.
$$
Finally, the last integral on the r.h.s. of (6.33) can be estimated by
$$\eqal{
&c\bigg[\intop_{\Omega_\varepsilon}\bigg|{\psi_{,r}\over r}
{\zeta_{,r}^2\over r}\bigg|^2r^2dx+\intop_{\Omega_\varepsilon}
\bigg|{\psi_{,z}\over r}{\zeta_{1,r}^2\over r}\bigg|^2r^2dx+
\intop_{\Omega_\varepsilon}\bigg|\psi\bigg({\zeta_{1,r}^2\over r^2}\bigg)_{,r}
\bigg|^2r^2dx\cr
&\le c(r_0)\bigg[\intop_\Omega(v_r^2+v_z^2)dx+\intop_\Omega\psi^2dx\bigg]\le
c(1/r_0)\intop_\Omega(v_r^2+v_z^2)dx.\cr}
$$
Therefore, (6.33) implies
$$
\intop_{\Omega_\varepsilon}|\nabla\tilde v_{z,r}|^2dx\le cI_1+c(1/r_0)
\intop_\Omega(v_r^2+v_z^2)dx.
\leqno(6.35)
$$
From (6.25), (6.26) and (6.35) after passing with $\varepsilon$ to 0 we 
have that
$$
\|\tilde v'\|_{L_2(0,T;H^2(\Omega))}\le cA_0,
\leqno(6.36)
$$
where
$$
A_0=\varphi\bigg(d_1,d_2,{1\over r_0}\bigg)\bigg(1+
\bigg\|{\tilde v_\varphi^2(0)\over r}\bigg\|_{L_2(\Omega)}+
\bigg\|{\tilde\chi(0)\over r}\bigg\|_{L_2(\Omega)}\bigg)
$$
Similarly, we have
$$
\|\tilde v'\|_{L_\infty(0,T;H^1(\Omega))}\le cA_0.
\leqno(6.37)
$$
This concludes the proof of Lemma 6.4.

Finally, we need

\proclaim Lemma 6.5. Assume that 
${\tilde\chi\over r}\in V_2^0(\Omega_\varepsilon^t)$, 
$v'\in V_2^0(\Omega_\varepsilon^t)$, $v'=(v_r,v_z)$. Then
$$
\|\tilde v/r\|_{V_2^1(\Omega_\varepsilon^t)}\le c
(\|\tilde\chi/r\|_{V_2^0(\Omega_\varepsilon^t)}+c(1/r_0)
\|v'\|_{V_2^0(\Omega_\varepsilon^t)}),
\leqno(6.38)
$$
where $t\le T$.

\Proof 
From (6.23) we have
$$
\intop_{\Omega_\varepsilon}\bigg|\nabla
\bigg({\tilde v_r\over r}\bigg)_{,r}\bigg|^2dx\le c
\intop_{\Omega_\varepsilon}(\vartheta^2+\vartheta_{,r}^2)dx.
$$
In view of definition of $\vartheta$ (see (4.6) and (4.7)) we derive
$$
\intop_{\Omega_\varepsilon}(\vartheta^2+\vartheta_{,r}^2)dx\le c(1/r_0)
\intop_{\Omega_\varepsilon}(v_z^2+v_{z,r}^2)dx+c\intop_{\Omega_\varepsilon}
\bigg(\bigg|{\tilde\chi\over r}\bigg|^2+
\bigg|\bigg({\tilde\chi\over r}\bigg)_{,r}\bigg|^2\bigg)dx.
$$
From (6.26) we obtain
$$
\intop_{\Omega_\varepsilon}{\tilde v_{r,zz}^2\over r^2}dx\le c
\intop_{\Omega_\varepsilon}\vartheta_{,z}^2dx
$$
where
$$
\intop_{\Omega_\varepsilon}\vartheta_{,z}^2dx\le\intop_{\Omega_\varepsilon}
\bigg({\tilde\chi\over r}\bigg)_{,z}^2dx+c(1/r_0)
\intop_{\Omega_\varepsilon}(v_r^2+v_{z,z}^2)dx.
$$
In view of the above inequalities we obtain
$$\eqal{
&\intop_{\Omega_\varepsilon}\bigg|\nabla^2
\bigg({\tilde v_r\over r}\bigg)\bigg|^2dx\le c\intop_{\Omega_\varepsilon}
\bigg(\bigg|{\tilde\chi\over r}\bigg|^2+
\bigg|\nabla{\tilde\chi\over r}\bigg|^2\bigg)dx\cr
&\quad+c(1/r_0)\intop_{\Omega_\varepsilon}(v_r^2+v_z^2+|\nabla v_z|^2)dx.\cr}
\leqno(6.39)
$$
Inequality (4.18) yields
$$
\intop_{\Omega_\varepsilon}|\nabla\eta_{,z}|^2dx\le
\intop_{\Omega_\varepsilon}\vartheta^2dx,
\leqno(6.40)
$$
where the l.h.s. of (6.40) equals
$$
\intop_{\Omega_\varepsilon}\bigg|\nabla{\psi_{,z}\over r^2}\bigg|^2dx=
\intop_{\Omega_\varepsilon}\bigg|\nabla{\tilde v_r\over r}\bigg|^2dx
$$
and the r.h.s. of (6.40) is estimated by
$$
\intop_{\Omega_\varepsilon}\vartheta^2dx\le\intop_{\Omega_\varepsilon}
\bigg|{\tilde\chi\over r}\bigg|^2dx+c(1/r_0)\intop_{\Omega_\varepsilon}
v_z^2dx,
$$
where (6.34) was used.

\noindent
Hence, (6.40) takes the form
$$
\intop_{\Omega_\varepsilon}\bigg|\nabla{\tilde v_r\over r}\bigg|^2dx\le c
\intop_{\Omega_\varepsilon}\bigg|{\tilde\chi\over r}\bigg|^2dx+c(1/r_0)
\intop_{\Omega_\varepsilon}v_z^2dx.
\leqno(6.41)
$$
Adding the integral with respect to time of (6.39) to $L_\infty$-norm with 
respect to time of (6.41) implies (6.38). This concludes the proof.

\section{References}

\item{BIN.} Besov, O. V.; Il'in, V. P.; Nikolskii, S. M.: Integral 
representation of functions and theorems of imbedding, Nauka, Moscow 1975 
(in Russian).

\item{CL.} Chae, D.; Lee, J.: On the regularity of the axisymmetric solutions 
of the Navier-Stokes equations, Math. Z. 239 (2002), 645--671.

\item{Ko.} Kochin, N. E.: Vectorial calculations and beginning of tensor
calculations, Akademia Nauk SSSR, Moscow 1951 (in Russian).

\item{L1.} Ladyzhenskaya, O. A.: The mathematical theory of viscous 
incompressible flow, Nauka, Moscow 1970 (in Russian).

\item{L2.} Ladyzhenskaya, O. A.: On unique solvability of three-dimensional 
Cauchy problem for the Navier-Stokes equations under the axial symmetry, 
Zap. Nauchn. Sem. LOMI 7 (1968), 155--177 (in Russian).

\item{L3.} Ladyzhenskaya, O. A.: Solutions "in the large" of the 
nonstationary boundary value problem for the Navier-Stokes system with two 
space variables, Comm. Pure. Appl. Math. 12 (1959), 427--433.

\item{LL.} Landau, L.; Lifshitz, E.: Hydrodynamics, Nauka, Moscow 1986 
(in Russian).

\item{LSU.} Ladyzhenskaya, O. A.; Solonnikov, V. A.; Uraltseva, N. N.: 
Linear and quasilinear equations of parabolic type, Nauka, Moscow 1967 
(in Russian).

\item{LU.} Ladyzhenskaya, O. A.; Uraltseva, N. N.: Linear and quasilinear 
equations of elliptic type, Nauka, Moscow 1973 (in Russian).

\item{MTL.} Mahalov, A.; Titi, E. S.; Leibovich, S.: Invariant helical 
subspaces for the Navier-Stokes equations, Arch. Rat. Mech. Anal. 112 
(1990), 193--222.

\item{NZ.} Nowakowski, B.; Zaj\c aczkowski, W. M.: Global existence of 
solutions to Navier-Stokes equations in cylindrical domains, Appl. Math. 
(Warsaw) 36 (2) (2009), 169--182.

\item{RZ.} Renc\l awowicz, J.; Zaj\c aczkowski, W. M.: Large time regular
solutions to the Navier-Stokes equations in cylindrical domains, Top. Meth. 
Nonlin. Anal. 32 (2008), 69--87.

\item{SS.} Seregin, G.; \v Sver\'ak, V.: On type I singularities of the local 
axi-sym\-metric solutions of the Navier-Stokes equations

\item{SZ.} Seregin, G.; Zaj\c aczkowski, W. M.: A sufficient condition of 
regularity for axially symmetric solutions to the Navier-Stokes equations, 
SIAM J. Math. Anal. 39 (2007), 669--685.

\item{S.} Solonnikov, V. A.: Estimates of the solutions of a nonstationary 
linearized system of Navier-Stokes equations, Trudy Mat. Inst. Steklov 70 
(1964), 213--317; English transl. Amer. Math. Soc. Trans. Ser. 2, 65 (1967), 
51--137.

\item{UY.} Ukhovskij, M. R.; Yudovich, V. I.: Axially symmetric motions of 
ideal and viscous fluids filling all space, Prikl. Mat. Mekh. 32 (1968), 
59--69 (in Russian).

\item{ZZ.} Zadrzy\'nska, E.; Zaj\c aczkowski, W. M.: Nonstationary Stokes 
system in Sobolev spaces.

\item{Z1.} Zaj\c aczkowski, W. M.: Global special regular solutions to the 
Navier\--Stokes equations in a cylindrical domain under boundary slip 
condition, Gakuto Intern. Ser., Math. Sc. Appl. 21 (2004), 1--188.

\item{Z2.} Zaj\c aczkowski, W. M.: Global special regular solutions to the 
Navier\--Stokes equations in a cylindrical domain without the axis of 
symmetry, Top. Meth. Nonlin. Anal. 24 (2004), 69--105.

\item{Z3.} Zaj\c aczkowski, W. M.: Global axially symmetric solutions with 
large swirl to the Navier-Stokes equations, Top. Meth. Nonlin. Anal. 29 
(2007), 295--331.

\item{Z4.} Zaj\c aczkowski, W. M.: Global regular axially symmetric solutions 
to the Navier-Stokes equations in a periodic cylinder, to be published.

\item{Z5.} Zaj\c aczkowski, W. M.: The H\"older regularity of swirl, to be 
published.

\item{Z6.} Zaj\c aczkowski, W. M.: On global regular solutions to the 
Navier-Stokes equations in cylindrical domains, Top. Meth. Nonlin. Anal. 
37 (2011), 55--85.

\item{Z7.} Zaj\c aczkowski, W. M.: Long time existence of regular solutions
to Navier-Stokes equations in cylindrical domains under boundary slip 
conditions, Studia Math. 169 (2005), 243--285.

\item{Z8.} Zaj\c aczkowski, W. M.: Global special regular solutions to the 
Navier\--Stokes equations in axially symmetric domains under boundary slip 
conditions, Diss. Math. 432 (2005), pp. 138.

\item{Z9.} Zaj\c aczkowski, W. M.: Global regular solutions to the 
Navier-Stokes equations in a cylinder, Banach Center Publ. 74 (2006), 
235--255.

\item{Z10.} Zaj\c aczkowski, W. M.: Existence and regularity properties of 
solutions of some elliptic system in domains with edges, Diss. Math. 274 
(1988), 95 pp.

\bye